\documentclass[a4paper,11pt, leqno]{article}
\usepackage[latin1]{inputenc}
\usepackage{hyperref}
\usepackage{articlemacro}

\usepackage{relsize}
\usepackage[bbgreekl]{mathbbol}
\usepackage{amsfonts}

\DeclareSymbolFontAlphabet{\mathbb}{AMSb}
\DeclareSymbolFontAlphabet{\mathbbl}{bbold}
\newcommand{\Prism}{\mathbbl{\Delta}}

\usepackage{tikz}
\usetikzlibrary{arrows,positioning,decorations.markings}
\tikzset{>=angle 90}

\begin{document}

\title{On the work of Peter Scholze}

\author{T.~Wedhorn}

\maketitle


\noindent{\scshape Abstract.\ }

This is a survey article over some of the work of Peter Scholze for the Jahresbericht der DMV. No originality is claimed.


\section*{Introduction}

In August 2018, Scholze was awarded a Fields medal ``for transforming arithmetic algebraic geometry over $p$-adic fields through his introduction of perfectoid spaces, with application to Galois representations, and for the development of new cohomology theories''.

Scholze's theory of perfectoid spaces together with its applications has constructed a new bridge between characteristic zero and positive characteristic. This allowed him to prove a number of longstanding conjectures in Algebraic Geometry and Number Theory. Examples are a conjecture of Ash within the Langlands program on Galois representations attached to torsion cohomology classes of locally symmetric spaces and Deligne's weight monodromy conjecture for complete intersections. Together with Bhatt, he developed the formalism of prismatic cohomology that includes an integral version of $p$-adic Hodge theory providing new links between different cohomology theories.

By now the richness and audacity of the theoretic edifice built on Scholze's ideas by Scholze himself and by many others might have become somewhat intimidating for beginners in its vastness. This text intends to give a guideline to many of these fascinating ideas.

It starts with an introductory part that gives several examples why geometry and arithmetic in characteristic $p$ -- beyond its inherent fascination -- is essential in modern theoretical mathematics. Moreover classical pre-Scholze (and still very important) ways to bridge the gap between characteristic zero and characteristic $p$ are described.

In the second part the new bridge of perfectoid spaces and diamonds is explained. It also includes two important direct applications (the weight monodromy conjecture for complete intersections and Hochster's conjecture) and a section on a central tool of the theory, the pro-\'etale topology.

The third part surveys some results of Scholze (often with coauthors) on $p$-adic Hodge theory. It also contains a longer section on prismatic cohomology, developed by Bhatt and Scholze, that generalizes several $p$-adic cohomology theory and goes far beyond the scope of $p$-adic Hodge theory.

In the fourth part we review some of Scholze's result within the net of theorems and conjectures that is called the Langlands program. It starts with a very brief introduction to the Langlands correspondence and two of its major players, moduli spaces of shtukas and Shimura varieties. Needless to say that those few pages do not contain any precise statements but only hints how one might think about these objects.

The last part lists two further achievements of Scholze and his coauthors. The first one is foundational work on topological cyclic homology by Nikolaus and Scholze. The second, very recent achievement together with Clausen is nothing less than a new definition what a (topological) space should be, introducing the new notion of condensed objects for an arbitrary ($\infty$-)category yielding an improved categorical framework for topological algebra.

This survey is far from complete. Although the text became much longer than initially planned, it contains still quite a few serious gaps. This is hardly surprising considering that the mathematical work of Scholze with various coauthors since the appearance of his first preprint in 2010 now surpasses by far 2000 pages of densely written mathematics.

For instance, I neglected almost totally the fundamental work of Bhatt and Scholze on the Witt vector affine Grassmannian (\cite{BS_WittGrass}), Scholze's ideas about the $q$-deformation of de Rham cohomology (\cite{Scholze_qDef} and \cite{BS_Prismatic}), or many of the fascinating results on period rings in \cite{Scholze_PAdicHodge}. Moreover, some topics are touched only superficially due to lack of space or competence of the author (e.g., Section~\ref{TCH}). The educated reader will certainly find many more serious omissions.

This is not a historical account. I did not try to explain any of the previous work on which Scholze's mathematics is built. I also ignored very often contributions of other mathematicians to the theory. In particular, I am aware of the fact that I did not give sufficient credit in many places. I tried to give precise references as often as possible but very often only to a most recent version of a result that hides the contributions of many mathematicians. Please refer to the references given in my references for this. I also took the risk to include many ``philosophical'' remarks about why one considers certain objects, what are the principal properties, or what are the main obstacles. Of course, all these remarks should be viewed as personal opinion of the author and not as definitive wisdom -- or wisdom at all.

\bigskip

\textsc{Acknowledgement}: I am very grateful to T.~Henkel, J.~Hesse, T.~Richarz and C.~Yaylali for many helpful remarks. Special thanks go to M.~Rapoport for sending me his laudation for Scholze's Fields medal \cite{Rapoport_Scholze}, which was a very helpful guideline, and for his comments on a first version of the text. Finally, I am grateful for the remarks of an anonymous referee.

\subsection*{Prerequisites}

Besides basic knowledge in Algebra, Analysis, Category Theory, Geometry, and Topology some knowledge of Algebraic Geometry will be assumed. For an overview of the basic notions in Algebraic Geometry I recommend \cite{Gortz_EGA}. It is hoped that Part~I and~II are readable for a broader audience at least if one accepts that a basic principle of geometry is to attach to geometric spaces objects in linear algebra via (co)homology theories. For many of the more arithmetic sections in particular in Part IV also some basic number theory in the sense of Weil \cite{Weil_BasicNT} will be very helpful.


\newpage

\tableofcontents


\subsection*{Notation}

The letter $p$ always denotes a prime number. We write $\FF_p$ for the field with $p$ elements.

All rings are assumed to be commutative. The group of units of a ring $R$ is denoted by $R^{\times}$. A ring $R$ is called \emph{local} if it has a unique maximal ideal $\mfr$. In this case $R/\mfr$ is called the \emph{residue field} of $R$. More generally, if $R$ is an arbitrary ring and $\pfr \subset R$ is a prime ideal, then the localization $R_{\pfr}$ of $R$ at $\pfr$ is a local ring with maximal ideal $\pfr R_{\pfr}$ and its residue field is called the \emph{residue field at $\pfr$}.

The topological groups that we will encounter will be usually abelian (most of the time they will be underlying additive groups of a ring) and first countable. In such topological groups we say that a sequence $(x_n)_n$ is a \emph{Cauchy sequence} if for every neighborhood $U$ of $0$ there exists an $N \in \NN$ such that $x_n - x_m \in U$ for all $n,m \geq N$. Such a topological groups is called \emph{complete} if it is Hausdorff and if every Cauchy sequence converges.

We use ``variety'' as a generic and imprecise way of a scheme that has sufficiently good properties. Varieties are always assumed to be separated and of finite type over their base but often have additional properties such as being reduced, quasi-projective, or geometrically connected. Sometimes we indicate the precise definition in a footnote.

By limits always projective or inverse limits are meant. Inductive limits are called colimits.

\subsubsection*{Global and local fields}

The fields of most interest for algebraic number theory are global and local fields. Both notions can be defined in characteristic zero and in positive characteristic. In characteristic zero, a global field is simply a finite extension of the field of rational numbers $\QQ$ and a local field is either equal to $\RR$ or to $\CC$ (archimedean case) or to a finite extension of the field of $p$-adic numbers $\QQ_p$ (non-archimedean case). In characteristic $p$, a global field is a finite extension of $\FF_p(t)$, the field of fractions of the polynomial ring $\FF_p[t]$, and a local field is a finite extension of $\FF_p((t))$, the field of fractions of the ring $\FF_p[[t]]$ of formal power series.


\section*{Some Background}

It is strongly suggested to skip this section and only refer to it when necessary.

\subsection*{Grothendieck topologies}

To define a sheaf on a topological space $X$ one needs only the category of open subsets of $X$ to define a presheaf to be a functor from $\textrm{Open}(X)^{\opp}$ to the category of sets (or of abelian groups, or...) and the notion of an open covering of an open subset to formulate the condition for a presheaf to be a sheaf. This can be formalized by defining a so-called \emph{Grothendieck topology} on a category by specifying certain families of morphisms in this category as coverings. Intersections of open sets are replaced by fiber products that are assumed to exist. Categories endowed with a Grothendieck topology are called \emph{sites}\footnote{We follow here the convention of the Stacks project \cite{Stacks} which differs from the original definition in \cite{SGA4}.}. We refer to \cite[Part 1, 2.3]{FGAExplained} for a nice introduction.

In Algebraic Geometry the category is usually the category of certain schemes and the topology is generated\footnote{``Generated'' means that one has to check the sheaf condition only for these special coverings.} by Zariski covering $(U_i \to U)_i$, i.e., $U_i \subseteq U$ are open subschemes with $\bigcup U_i = U$, and coverings consisting of a single surjective morphism $f\colon \Spec B \to \Spec A$ of affine schemes having a certain fixed property $\Pbf$. For instance, if $\Pbf$ is the property ``isomorphism'' (resp.~``\'etale'', resp.~``flat''), one obtains the Zariski topology (resp.~the \'etale topology, resp. the fpqc-topology)\footnote{Here we ignore some set-theoretic issues.}. For instance, the (small) \'etale site of a scheme $S$ is the category of all $S$-schemes $f\colon U \to S$ with $f$ \'etale endowed with the topology given by covering $(f_i\colon U_i \to U)_{i\in I}$, where $f_i$ is \'etale and $\bigcup f_i(U_i) = U$.

\subsection*{Cohomology of sheaves}

Once a Grothendieck topology is fixed one has a formalism of sheaves and their cohomology. See the \cite[Part 1, Chap.~2]{FGAExplained} for a brief introduction of the notion of a sheaf on Grothendieck site and \cite{KS_CategoriesSheaves} or \cite{Stacks} for a more encyclopedic treatment. The cohomology of an abelian sheaf is defined as the derived functor of the functor of global sections.

\subsection*{Reductive groups}

In the fourth part about the Langlands program we will also use the notion of a reductive group. Let us briefly recall what this means. A group object in the category of varieties over a field $k$ is called an \emph{algebraic group} over $k$. Such an algebraic group is called \emph{linear} if there exists an embedding of algebraic groups into $\GL_n$ for some $n$. If $k$ is algebraically closed, then a linear algebraic group over $k$ is called \emph{reductive} if it is connected and if the largest normal subgroup scheme of $G$ that can be embedded into a group of unipotent matrices is trivial. A linear algebraic group over an arbitrary field $k$ is called \emph{reductive} if its base change to some algebraically closed extension of $k$ is reductive.

Every linear algebraic group over a field of characteristic zero can be defined over a subfield of $\CC$. By definition these groups are reductive if and only if their base change to $\CC$ or $\RR$ is reductive. A connected complex linear algebraic group $G$ is reductive if and only if it has a compact real form, i.e., there exists an automorphism $\theta$ of $G$ of order $2$ such that $\set{g \in G(\CC)}{\theta(\gbar) = g}$ is a compact subgroup. A connected real linear algebraic group $G$ is reductive if and only if there exists an embedding $G \mono \GL_{n,\RR}$ such that $G$ is stable under transposition.

For instance, the groups $\GL_n$, $\Sp_n$, $\SO_n$ are reductive (in arbitrary characteristic). The group of upper triangular matrices is not reductive.

A linear algebraic group over a field of characteristic $0$ is reductive if and only if every finite-dimensional algebraic representation is a direct sum of simple representations.

A subgroup $P$ of a reductive group $G$ is called \emph{parabolic} if the quotient variety $G/P$ is projective. A parabolic subgroup $B$ is called \emph{Borel subgroup} if after base change to some algebraic closure $\kbar$, $B_{\kbar}$ is a minimal parabolic subgroup of $G_{\kbar}$. Borel subgroups always exist if $k$ is algebraically closed but not necessarily over general fields. If they exist, the reductive group is called \emph{quasi-split}.

Reductive groups over separably closed fields $k$ are essentially combinatorial data. To every reductive group over $k$ one attaches a based root datum consisting of a finitely generated free abelian group $X^*$, finite subsets $R$ of $X^*$, the set of roots, and $R\vdual$ of the dual $X_*$ of $X^*$, the set of coroots, and a system $\Delta \subseteq R$ of so-called simple roots. These data have to satisfy certain axioms which can be formulated with the language of Linear Algebra. The main classification theorem of reductive groups then states that this construction yields a bijection between isomorphism classes of reductive groups over $k$ and based root data. This in particular implies that over an algebraically closed field a linear algebraic group $G$ is reductive if and only if its center $C$ can be embedded into a group of diagonal matrices and $G/C$ is isomorphic to a product of $\PGL_n$ ($n \geq 1$), ${\rm PSO}_n$ ($n \geq 6$), ${\rm PSp}_{2n}$ ($n \geq 2$), or of one of five exceptional groups.

If $G$ is a reductive group over a general field $k$ with separable closure $\kbar$, then the root datum of $G_{\kbar}$ carries a continuous action of $\Gal(\kbar/k)$ that preserves the set of simple roots if $G$ is quasi-split. This yields a bijection between isomorphism classes of quasi-split reductive groups over $k$ and isomorphism classes of based root data with continuous $\Gal(\kbar/k)$-action.


\part{A wonderful positive world}

\section{Positive characteristic}

Recall that $p$ always denotes a prime number. The geometry and arithmetic of varieties over fields of characteristic $p$ has often a very different flavor than of those over the complex numbers. In fact, in many respects geometric and arithmetic problems are easier to handle in characteristic $p$ than in characteristic zero. The main reason is the existence of an additional symmetry, namely the Frobenius endomorphism $x \sends x^p$. It is for rings $R$ of characteristic $p$ (i.e. in which $p = 0$ holds) not only compatible with multiplication, i.e. $(ab)^p = a^pb^p$ for $a,b \in R$, but also with addition, i.e., $(a+b)^p = a^p + b^p$. In other words, it is a ring endomorphism. In a certain sense one can view it as a ``contracting'' endomorphism which is ubiquitous in characteristic $p$. 

\subsection{The Frobenius}


A non-zero ring $R$ is said to have \emph{characteristic $p$} if $p = 1 + 1 + \dots + 1$ ($p$-times) is zero in $R$. Equivalently, the unique ring homomorphism $\ZZ \to R$ factors through $\FF_p = \ZZ/p\ZZ$.

If $R$ is a ring of characteristic $p$, then the map
\[
\sigma\colon R \to R, \qquad x \sends x^p
\]
is a ring endomorphism and it is called the \emph{(arithmetic) Frobenius} of $R$. Sometimes we will call the map $x \sends x^q$ for some power $q$ of $p$ the $q$-Frobenius.

A ring of characteristic $p$ is said to be \emph{perfect} if the Frobenius map $x \sends x^p$ is bijective. Examples of perfect rings are finite fields and algebraically closed fields. Below we will also see examples of perfect rings that are not fields.

An example of a field that is not perfect is $k(T)$, the field of fractions of the polynomial ring $k[T]$ over a field $k$ of characteristic $p$. Here the indeterminate $T$ is not the $p$-th power of any polynomial.

Let us now give some examples for the usefulness of the Frobenius symmetry. They range from using methods of characteristic $p$ only as a tool -- though an essential one -- to obtain results over the complex numbers in Example~\ref{MORI} over arithmetic questions that are interesting over any field but often need methods of positive characteristic to get off the ground in Examples~\ref{CLASSQUAD} and~\ref{COUNTPOINT} to the Langlands program that is really an entire web of theorems, conjectures, and principles linking large parts of theoretical mathematics.

\subsection{Counting Points}\label{COUNTPOINT}

Systems of polynomial equations in several variables and finding their solution is an omnipresent problem in and outside mathematics and has been studied since antiquity. In modern language that means, given a variety $X$ over a field $k$, we are interested in determining the set $X(k)$ of $k$-rational points. For instance, if $X$ is given as a vanishing set of polynomials $f_1,\dots,f_r$ in $n$ variables with coefficients in $k$, then $X(k)$ is the set of solutions $x \in k^n$ of the system of equations $f_1(x) = \dots = f_r(x) = 0$.

To determine $X(k)$ is in general a very hard question as already the ``simple'' case $k = \QQ$, with one equation in two variables $f = X^m + Y^m - 1$ for $m \geq 3$ shows. The question, whether $f$ has a solution $(x,y) \in \QQ^2$ with $x,y \ne 0$ is equivalent to Fermat's last theorem by clearing denominators.

Things get much easier over finite fields due to Grothendieck-Lefschetz theory. Let $k = \FF_q$ be a field with $q$ elements, where $q$ is a power of $p$, and let $X$ be a variety\footnote{a separated scheme of finite type} over $\FF_q$. Then $X(\FF_q)$ is a finite set and one has
\begin{equation}\label{EqGrothLefschetz}
\#X(\FF_q) = \sum_{i=0}^{2\dim(X)}(-1)^i\tr(F_q\mid H^i_c(X)).
\end{equation}
Here $F_q$ is the geometric $q$-Frobenius\footnote{If $\kbar$ is an algebraic closure of $k$ and $X(\kbar) = \set{x \in \kbar^n}{f_1(x) = \dots = f_r(x) = 0}$ for polynomials $f_i \in \FF_q[T_1,\dots,T_n]$, then $F_q(x_1,\dots,x_n) = (x_1^q,\dots,x_n^q)$ for $x = (x_1,\dots,x_n) \in \kbar^n$.} acting as a linear endomorphism on the cohomology\footnote{$\ell$-adic cohomology with compact support for some auxiliary prime number $\ell \ne p$} of the base change of $X$ to some algebraic closure of $k$. This yields a cohomolo\-gical way to calculate the number of rational points which has turned out to be very useful.

For instance, a striking application was given by B.~C.~Ng\^o in his proof \cite{Ngo_Fondamental} of the Fundamental Lemma for which he was awarded a Fields medal in 2010. His proof translated a purely combinatorial statement (the Fundamental Lemma) into a statement on counting points over a finite field and then used \eqref{EqGrothLefschetz} to conclude\footnote{Needless to say that this description of the proof sweeps a whole world of technical complications under the rug. Note that by now there is also a quicker proof of the main technical point of the Fundamental Lemma using $p$-adic integration by Groechenig, Wyss, and Ziegler \cite{GWZ}.}.


\subsection{Classification of quadratic spaces}\label{CLASSQUAD}

To classify quadratic spaces (i.e., finite-dimensional vector spaces endowed with a quadratic form) over the field $\QQ$ of rational numbers, one proceeds usually as follows (\cite{Serre_CourseArithmetic}). By the theorem of Hasse-Minkowski, two quadratic spaces over $\QQ$ are isomorphic if and only if they are isomorphic over $\QQ_v$ for every place $v$ of $\QQ$. Here $\QQ_v$ denotes the completion of $\QQ$ with respect to the metric given by a non-trivial absolute value $v$ of $\QQ$. For the usual archimedean absolute value $v$ we obtain $\QQ_v = \RR$ and quadratic spaces over the real numbers are classified by dimension and signature by Sylvester's law of inertia. 

The other places $v$ of $\QQ$ are (up to equivalence) the $p$-adic absolute values $|\cdot|_p$ for each prime number $p$. For $|\cdot|_p$ a number is the closer to $0$ the larger a power of $p$ can be extracted. The completion of $\QQ$ with respect to the $p$-adic absolute value is the field $\QQ_p$ of $p$-adic numbers. As the $p$-adic absolute value is non-archimedean (i.e., $|a+b|_p \leq \max\{|a|_p,|b|_p\}$ for $a,b \in \QQ_p$), one has the subring $\ZZ_p = \set{a \in \QQ_p}{|a|_p \leq 1}$ of $p$-adic integers which is compact. It can be also written as projective limit
\[
\ZZ_p = \lim_n \ZZ/p^n\ZZ.
\]
The ``open unit ball'' $\set{a \in \QQ_p}{|a|_p < 1}$ is the unique maximal ideal of $\ZZ_p$. It is the ideal generated by $p$. The residue field $\ZZ_p/p\ZZ_p$ is $\FF_p$, the field with $p$ elements.

Now quadratic spaces over $\QQ_p$ can be completely classified by dimension, discriminant, and Hasse invariant. Moreover, suppose the quadratic form is given by 
\[
q(x) = a_1x_1^2 + \dots + a_rx_r^2 + pa_{r+1}x_{r+1}^2 + \dots + pa_nx_n^2
\]
for units $a_1,\dots,a_r,\dots,a_n \in \ZZ_p^{\times}$ (this can always be achieved by some linear base change). Let $\abar_i$ the image of $a_i$ in $\FF_p^{\times}$. Then for $p \ne 2$ the two quadratic forms
\[
\qbar_1(y) = \abar_1y_1^2 + \dots + \abar_ry_r^2, \qquad \qbar_2(z) = \abar_{r+1}z_{1}^2 + \dots + \abar_nz_{n-r}^2
\]
over $\FF_p$ determine $q$ uniquely up to isometry.

Hence in the end, the classification of rational quadratic spaces is based on the classification of quadratic spaces over finite fields. But this is quite easy: Two non degenerate quadratic spaces over a finite field are isomorphic if and only if they have the same dimension and discriminant.

\subsection{Producing rational curves in the Minimal Model Program}\label{MORI}

The Minimal Model Program strives to classify algebraic varieties over the complex numbers and tackles therefore one of the central questions in modern Algebraic Geometry. Its goal is to construct for every smooth projective variety $X$ over the complex numbers a birational model that is ``as simple as possible''. For instance, the line bundle of volume forms should be sufficiently positive, more precisely, the canonical divisor should be nef. For this it is important to identify simple obstructions that prevent the positivity. These are rational curves, i.e. non-constant morphisms $\PP^1 \to X$, of sufficiently small degree, on which the canonical divisor is negative. One would like to get rid of them by contracting them.

Essentially by definition of not being nef there always exists a curve $\gamma\colon C \to X$ on which the canonical divisor is negative. Now the idea is to use ``Bend \& Break'', that is to deform (``bend'') the curve $C$ within $X$ until it degenerates into several curves of smaller genus and degree (it ``breaks'').  A calculation using the scheme of all morphisms $C \to X$ yields that such a deformation is possible if the pullback $\gamma^*T_X$ of the tangent bundle of $X$ to $C$ has sufficiently negative degree. This can for instance be achieved by precomposing $\gamma$ with a map of curves of high degree without changing the genus of the curve too much.

Over the complex numbers it is not clear how to find such maps and there is no proof of this using only methods in characteristic $0$. But in characteristic $p$ one can simply precompose $\gamma$ with sufficiently high powers of the Frobenius endomorphism of $C$. Hence one uses a standard technique to find models of $X$ and $C$ in characteristic $p$ for all but finitely many prime numbers $p$ (see also Subsection~\ref{ZEROTOP} below), then proves a Bend \& Break lemma there using the Frobenius, and finally uses the fact that the lemma holds in characteristic $p$ for infinitely many prime numbers $p$ to conclude that it also holds over the complex numbers (see~\cite[1.10]{KollarMori} for more details). This gives an example where the Frobenius is an indispensable tool, even if one is interested only in questions in characteristic zero. There are many more such examples.

\subsection{The Langlands program}\label{INTROLANGLANDS}

The Langlands program is an entire web of theorems, conjectures, and principles about connections between number theory, geometry, representation theory, and analysis. It relates Galois groups over local and global fields (see the introduction for our use of the notion of local and global fields) to representations of Lie groups over these fields. 

In its most basic form, the Langlands conjecture predicts -- roughly speaking -- for all $n \geq 1$ a correspondence between certain $n$-dimensional representations of the Galois group of a global field $F$ (resp.~a local field $K$) and certain representations of $\GL_n(\AA_F)$ (resp.~of $\GL_n(K)$). Here $\AA_F$ denotes the ring of adeles of the global field $F$. If $n = 1$, then the Langlands correspondence is equivalent to class field theory of the global or local field. There exist also generalizing conjectures, where $\GL_n$ is replaced by more general algebraic groups.

For global fields of characteristic zero, the core of number theory, the Langlands conjecture is still wide open. In all other cases there exist now theorems. Again the case of characteristic $p$ has turned out to be more accessible essentially because of two reasons. The first one is that global fields in characteristic $p$ are the same as function fields of curves over finite fields and this allows to apply more geometric methods to the Langlands correspondence. The other reason is again the existence of the Frobenius which is an indispensable tool in all proofs of the Langlands correspondence in positive characteristic.

Below in Subsection~\ref{GEOMLLC} we will see how the ideas of Scholze also allow to carry over these ideas to a quite general form of the Langlands conjecture for non-archimedean local fields of characteristic zero.


\section{The classical way to pass between zero and positive characteristic}

\subsection{Passage from characteristic zero to positive characteristic}\label{ZEROTOP}

\subsubsection{Integral models over $\ZZ$}

Suppose that there is given a variety $X$ over the field of rational numbers $\QQ$ (or over some finite extension of $\QQ$). For simplicity let us assume that $X$ is projective, i.e., $X$ is given as the vanishing locus of homogeneous polynomials $f_1,\dots,f_r$ with rational coefficients in a projective space
\[
X = \set{x = (x_0 : \dots : x_n) \in \PP^n_{\QQ}}{f_1(x) = \dots = f_r(x) = 0}.
\]
As each of these polynomials has only finitely many coefficients, we may multiply each of the polynomials with some common denominator to obtain polynomials with integer coefficients. This does not change $X$ and hence we may assume that all $f_i \in \ZZ[t_0,\dots,t_n]$. 

Let $p$ be a prime number. Then we can consider for each $f_i$ its image $\fbar_i$ in $\FF_p[t_0,\dots,t_n]$. The vanishing locus $\Xscr_p$ of $\fbar_1,\dots,\fbar_r$ in the projective space over $\FF_p$ is then a projective variety over $\FF_p$. In the language of modern Algebraic Geometry we consider the vanishing locus $\Xscr$ of $f_1,\dots,f_r$ in the projective space $\PP^n_{\ZZ}$ over the ring $\ZZ$ as a projective scheme over $\ZZ$. Then $\Xscr_p$ is simply the fiber of $\Xscr$ over the point $(p) \in \Spec \ZZ$.

Unfortunately, the model $\Xscr$ and its fiber $\Xscr_p$ over $\FF_p$ depend very much of the choice of the equations $f_1,\dots,f_r$ and in general there is no single optimal choice. Moreover, sometimes good properties of $X$ cannot be transferred to $\Xscr_p$. For instance, if $X$ is smooth over $\QQ$ (i.e., it has no singularities, or equivalently, if viewed as a complex analytic subspace of $\PP^n_{\CC}$ given by the same equations, then this complex analytic space is a manifold), it might not be possible for a given $p$ to find a model $\Xscr$ such that the $\FF_p$-variety $\Xscr_p$ is smooth.

\subsubsection{$p$-adic topology and $p$-adic numbers}\label{PADIC}

To study such phenomena in more detail one usually passes to $p$-adic numbers. For this let us recall the $p$-adic topology on a ring as it will also play an important role later.

We consider a somewhat more general notion. Let $R$ be a ring and let $I \subseteq R$ be an ideal. For $n \geq 0$ we denote by $I^n$ its $n$-th power, i.e., the ideal generated by all $n$-fold products of elements of $I$. We define a pseudo-metric\footnote{A pseudo-metric $d$ on a set $X$ is a distance function $d\colon X \times X \to \RR_{\geq 0}$ that satisfies all properties of a metric except the implication $d(x,y) = 0 \implies x = y$.} on $R$ by
\[
d(a,b) := \alpha^{\delta(a-b)},\qquad\text{with}\quad \delta(a) := \sup\set{n \in \NN_0}{a \in I^n}
\]
where $\alpha$ is some real number with $0< \alpha < 1$. The pseudo-metric induces a first countable topology on $R$ that makes addition and multiplication continuous and which is independent of the choice of $\alpha$. A basis of neighborhoods of $0$ is given by the family $(I^n)_{n\geq1}$. It is called the \emph{$I$-adic topology}. It is Hausdorff if and only if $d$ is a metric, i.e., if and only if $\bigcap_n I^n = 0$. If $I$ is a principal ideal generated by one element $\pi$ (for instance $\pi = p$), we call the $I$-adic topology also the $\pi$-adic topology. In this case one has $\delta(a) = n$ if $a$ is divisible by $\pi^n$ but not by $\pi^{n+1}$. Hence an element of the ring is close to zero with respect to the $\pi$-adic topology if it is divisible by a high power of $\pi$.

One can also consider the completion of $R$ with respect to this pseudo-metric. It depends only on the $I$-adic topology and one obtains a complete metrizable topological ring $\Rhat$ together with a map $R \to \Rhat$ of topological rings. If $R$ is noetherian, then the topology on $\Rhat$ is the $\Ihat$-adic topology, where $\Ihat = I\Rhat$ is the ideal in $\Rhat$ generated by the image of $I$\footnote{Beware that in general the topology on $\Rhat$ is not the $I\Rhat$-adic topology and that $\Rhat$ might not even by complete with respect to the $I\Rhat$-adic topology.}.

If we now specialize to $R = \ZZ$ and $I = (p)$, the $p$-adic completion of $\ZZ$ is the ring $\ZZ_p$ of $p$-adic integers. It is a discrete valuation ring, i.e. a principal ideal domain with a unique (up to multiplication with units) prime element. In fact, $p$ is a prime element and hence the ideal generated by $p$ is the unique maximal ideal of $\ZZ_p$. The residue field $\ZZ_p/(p)$ is the field $\FF_p$ with $p$ elements. The field of fractions of $\ZZ_p$ is the field $\QQ_p$ of $p$-adic numbers. It carries the $p$-adic absolute value $|\cdot|_p$ defined as follows. As $\QQ_p$ is the field of fractions of a factorial domain with the single (up to units) prime element $p$, we can write every non-zero element of $\QQ_p$ as $up^n$ for $u \in \ZZ_p^{\times}$ a unit in $\ZZ_p$ and $n \in \ZZ$. Then $|up^n|_p := p^{-n}$ defines the $p$-adic absolute value\footnote{One could also define the $p$-adic absolute value by $|up^n|_p := \alpha^{n}$ for any $0 < \alpha < 1$ and get the same topology. The reason for choosing $\alpha = p^{-1}$ becomes more apparent if one studies global questions involving all $p$ simultaneously.}. The $p$-adic absolute value satisfies all properties of a \emph{non-archimedean absolute value} $|\cdot|$, i.e.
\begin{assertionlist}
\item
$|a| = 0$ if and only if $a = 0$,
\item
$|ab| = |a||b|$ for all $a,b$,
\item
$|a+b| \leq \max{|a|,|b|}$ for all $a,b$.
\end{assertionlist}
It is \emph{non-trivial} in the sense that there are non-zero elements $a$ with $|a| \ne 1$.

The big advantage of the ring $\ZZ_p$ compared to $\ZZ$ is the property that it is $p$-adically complete and one has $\ZZ_p = \lim_n \ZZ/p^n\ZZ$. This allows in many cases to approximate solutions of polynomial equations of $\ZZ_p$ by such solutions over the rings $\ZZ/p^n\ZZ$. For instance Hensel's lemma tells us that if a polynomial $f \in \ZZ_p[t]$ in one variable has a simple root over $\FF_p$, then it has already a root over $\ZZ_p$. Geometrically, this means that one should view $\Spec \ZZ_p$ as a very small thickening of $\Spec \FF_p$.

Instead of studying models of varieties $X$ over $\QQ$ one often considers them locally in $p$, i.e., one considers $X$ as a variety over $\QQ_p$ which is a field extension of $\QQ$ and studies the question whether they have ``good'' integral models $\Xscr$ over the complete discrete valuation ring $\ZZ_p$. The fiber of such an integral model in characteristic $p$ is then a variety over $\FF_p$.

\subsubsection{Bridges towards characteristic $p$: Non-archimedean fields and their rings of integers}\label{NONARCH}

As it is important to study for many number theoretic questions equations over extensions of $\QQ$ one also considers more generally extensions $K$ of $\QQ_p$ and models over integral domains $R$ with $K = \Frac R$ and with a unique maximal ideal $\mfr$ such that $R/\mfr$ is a field of characteristic $p$. Such rings $R$ are examples of rings of \emph{mixed characteristic} (i.e., they contain prime ideals whose residue fields have different characteristics), in this case only characteristic $0$ and characteristic $p$ for a single prime number $p$.

These rings form a ``bridge'' between characteristic zero and characteristic $p$, so to speak. As to keep the ``distance between characteristic zero and characteristic $p$ as small as possible'' one usually assumes that the ring is $p$-adically complete. Among these rings $R$ the ``ultimate bridge'' is $\ZZ_p$ in the sense that there always exists a unique map of topological rings $\ZZ_p \to R$.

The standard way to use these bridges is the following. Let $K$ be a \emph{non-archimedean field}, i.e., a complete topological field whose topology is induced by a non-trivial non-archimedean absolute value $|\cdot|$ (\ref{PADIC}). As the absolute value is non-archimedean, the subset $O_K := \set{a \in K}{|a| \leq 1}$ is a subring of $K$ whose field of fractions is $K$. It has a unique maximal ideal, namely $\mfr = \set{a \in K}{|a| < 1}$. Moreover, $O_K$ is a \emph{valuation ring}, i.e., for every element $0\ne x \in K$ one has $x \in O_K$ or $x^{-1} \in O_K$. The topology on $O_K$, induced by the topology on $K$, makes $O_K$ into a complete topological ring endowed with the $\varpi$-adic topology, where $\varpi$ is any non-zero element in $\mfr$.

Now suppose that $K$ is an extension of $\QQ_p$ and that $|\cdot|$ extends the $p$-adic absolute value $|\cdot|_p$ on $\QQ_p$. For instance, for every algebraic extension $K$ of $\QQ_p$ there exists a unique extension of the $p$-adic absolute value which is given on a finite subextension $M$ by $|a| = (|N_{M/\QQ_p}(a)|_p)^{1/[M:\QQ_p]}$, where $N_{M/\QQ_p}\colon M \to \QQ_p$ is the norm. As $|p|_p = p^{-1} < 1$, one has $p \in \mfr$ and hence the residue field $k := O_K/\mfr$ is a field of characteristic $p$. Therefore $O_K$ is of mixed characteristic $(0,p)$ and $p$-adically complete.

Now that we have the ``bridge'' $R$, we can use it to reduce varieties in characteristic zero (non-canonically) to characteristic $p$. Let $X$ be a variety over $K$ that we assume to be projective for simplicity. Such varieties may for instance come by extension of scalars from a variety over a number field, because every number field can be embedded into a finite extension of $\QQ_p$. Now find a model $\Xscr$ of $X$ over $O_K$, i.e., a scheme $\Xscr$ over $O_K$ such that $\Xscr \otimes_{O_K} K \cong X$. Then we can reduce $\Xscr$ modulo $\mfr$ to obtain a scheme $\Xscr_k$ over $k$ which is now a scheme in characteristic $p$.

One can hope that algebraic invariants of $X$ and $\Xscr_k$ are closely related. For this one needs $\Xscr$ to be ``as nice as possible''. This includes that $\Xscr$ should be projective and flat over $O_K$ which is always possible. Moreover it should have singularities that are as mild as possible. This is a hard problem that by now has only a really satisfying solution if $X$ is a curve\footnote{Although using the language of log-schemes there has been much progress in higher dimension as well.}. For instance, if $X$ is smooth, we may hope that we can choose $\Xscr$ in such a way that $\Xscr_k$ is also smooth. If this is possible (which is not always the case), then we say that $X$ \emph{has good reduction}. In this case $X$ and $\Xscr_k$ share a number of common algebraic invariants, for instance they have the same $\ell$-adic cohomology (for a brief introduction to $\ell$-adic cohomology see Subsection~\ref{PROETSCHEME} below).


\subsection{Passage from positive characteristic to characteristic zero}\label{PTOZERO}

In the mathematical interplay between characteristic zero and characteristic $p$ one also would like to walk the other way. Thus one starts with objects in characteristic $p$ and wants to lift them to characteristic $0$. Hence one starts with a variety $\Xcal$ over a field $k$ of characteristic $p$. Then one first looks for a bridge towards characteristic zero, i.e. for a mixed characteristic local integral domain $A$ whose residue field is $k$ and whose field of fractions $K$ is a field of characteristic zero. The ring $A$ is by no means unique, although there is a ``minimal'' $A$ that we describe below if $k$ is perfect.

As a second step one tries to find a scheme $\Xscr$ over $A$ such that $\Xscr \otimes_A k \cong \Xcal$. For a given $A$ this is often not possible, and sometimes it is not possible regardless how one chooses $A$. Often this leads to very interesting deformation problem of which we cannot say more here. If one was able to construct $\Xscr$, the generic fiber $\Xscr \otimes_A K$ is a variety over a field of characteristic zero.

\subsubsection{The ring of Witt vectors}\label{WITT}

For perfect fields $k$ and more general for perfect rings, there is a canonical ``one para\-meter deformation'' towards characteristic zero given by the Witt ring\footnote{Here we mean the ring of $p$-typical Witt vectors. It also exists for non-perfect rings but has to be defined differently. We will not need this generality.}. It is a functor, that sends a perfect ring $R$ to a $p$-adically complete $p$-torsionfree $\ZZ_p$-algebra $W(R)$. We avoid the concrete somewhat lengthy construction of $W(R)$ (see \cite{Bou_AC} or \cite{Rabinov_Witt} for this), although it is important if one wants to work with Witt vectors. Instead, we will define $W(\cdot)$ as the left adjoint of another functor that is in some sense the fastest way to produce a perfect ring from a $p$-adically complete $\ZZ_p$-algebra. 

Hence let $B$ be in $\widehat{\textrm{(Alg)}}_{\ZZ_p}$, the category of $p$-adically complete $p$-torsionfree $\ZZ_p$-algebra. Then $B/pB$ is a ring of characteristic $p$. It is in general not perfect. One makes it perfect by setting
\begin{equation}\label{EqDefFlat}
B^{\flat} := \lim (B/pB \llefttoover{(\cdot)^p} B/pB \llefttoover{(\cdot)^p} B/pB\llefttoover{(\cdot)^p} \cdots).
\end{equation}
The functor $B \sends B^{\flat}$ from $\widehat{\textrm{(Alg)}}_{\ZZ_p}$ to the category of perfect rings has a left adjoint functor $R \sends W(R)$. The ring $W(R)$ is called the \emph{ring of Witt vectors of $R$}. The Frobenius automorphism on $R$ induces by functoriality an automorphism on $W(R)$, which is again called Frobenius. But beware, it is not given by $x \sends x^p$. One always has $W(R)/pW(R) = R$. For $B$ in $\widehat{\textrm{(Alg)}}_{\ZZ_p}$ the counit of the adjunction is a map
\begin{equation}\label{EqDefTheta}
\theta\colon W(B^{\flat}) \lto B.
\end{equation}
Moreover, the projection $W(R) \to R$ has a functorial multiplicative (but not additive) section $x \sends [x]$, and $[x] \in W(R)$ is called the \emph{Teichm\"uller representative} of $x$ in $R$.

If $k$ is a perfect field, then $W(k)$ is a complete discrete valuation ring with prime element $p$. The Witt ring $W(k)$ together with its isomorphism $W(k)/pW(k) \cong k$ is an initial object in the category of pairs $(R,\iota)$, where $R$ is a complete local noetherian ring and $\iota$ an isomorphism of its residue field with $k$. In this sense, $W(k)$ is the ``initial thickening of $k$ towards characteristic $0$''. For instance, one has $W(\FF_p) = \ZZ_p$.


\part{Perfectoid Spaces and Diamonds}

Scholze defined in his thesis the notion of a perfectoid ring $R$ and as a global variant the notion of a perfectoid space that are in general objects in mixed characteristic $(p,0)$. He also defined a tilting functor $X \sends X^{\flat}$ from the category of perfectoid spaces to the category of perfectoid spaces in characteristic $p$. Restricting the tilting functor to perfectoid spaces in characteristic zero gave a new and canonical way to pass from characteristic zero to characteristic $p$. Moreover, he showed that the categories of \'etale sheaves on $X$ and $X^{\flat}$ are equivalent, thus allowing to relate many geometric questions in characteristic $0$ and in characteristic $p$ directly.

In very special cases the tilting functor was already previously considered by Fontaine. For certain zero-dimensional spaces (i.e., certain fields) a theorem of Fontaine-Winten\-berger already proved the equivalence of the categories of \'etale sheaves over these points, which means for fields that the Galois groups of the fields are isomorphic. Scholze realized that these results were merely the simplest, zero-dimensional case of a general theory.

Of course, Scholze's construction of tilting does not make the ``classical way'' described above to pass between characteristic $0$ and $p$ obsolete. But it adds a fundamentally new and very useful tool. One ``problem'' with the theory of perfectoid spaces is the fact that almost no object classically considered in number theory is perfectoid. This does not make the theory useless because one can always ``perfectoidize''. Scholze showed that locally for a suitable topology (the pro-\'etale topology, see \ref{PROET} below) every classical $p$-adic variety is indeed perfectoid. In fact, starting with a $p$-adic variety, localizing to make it perfectoid, passing to its tilt, and then ``delocalizing'' this tilt can be made canonical. This is the theory of diamonds developed by Scholze.


\section{Scholze's way to pass between mixed and positive characteristic I: Perfectoid spaces}\label{PERFSPACES}

\subsection{Perfectoid rings}\label{PERFDRING}

The original notion of a perfectoid ring, introduced by Scholze in his thesis \cite{Scholze_Perfectoid}, has subsequently been generalized by Fontaine in \cite{Fontaine_BourbakiScholze} and Kedlaya and Liu in \cite{KedlayaLiu_I}. Further equivalent definitions have been also introduced in \cite{BMS_IntegralHodge} and \cite{BS_Prismatic}. By now there is a generally accepted definition of perfectoid rings recalled below.

Let $S$ be a ring that is $\pi$-adically complete for some $\pi \in S$ such that $\pi$ divides $p$. Then $S$ is also $p$-adically complete (\cite[090T]{Stacks}). As in \eqref{EqDefFlat} we define $S^{\flat}$, and we also set
\begin{equation}\label{EqDefAinf}
A_{\inf}(S) := W(S^{\flat}).
\end{equation}
Again there is a map $\theta\colon A_{\inf}(S) \to S$ as in \eqref{EqDefTheta}.

A ring $S$ is \emph{perfectoid} if it satisfies the following equivalent properties (see \cite[3.7--3.10]{BS_Prismatic}).
\begin{equivlist}
\item
The ring $S$ is $\pi$-adically complete for some $\pi \in S$ such that $\pi^p$ divides $p$, the Frobenius map $\phi\colon S/pS \to S/pS$ is surjective, and the kernel of $\theta\colon A_{\inf}(S) \to S$ is principal.
\item
The ring $S$ is isomorphic to $W(R)/(\xi)$, where $R$ is a perfect ring of characteristic $p$ and $\xi \in W(R)$ an element such that $\phi(\xi) - \xi^p$ is of the form $pu$ for some unit $u \in W(R)^{\times}$. Here $\phi$ denotes the Frobenius automorphism of $W(R)$.
\end{equivlist}
If these conditions are satisfied, then $R = S^{\flat}$ and the isomorphism $S \cong W(R)/(\xi)$ is given by $\theta$. For a perfectoid ring $S$, the perfect ring $S^{\flat}$ is called the \emph{tilt of $S$}.

The notion of a perfectoid ring and its tilt is also extended to \emph{complete Tate rings}, i.e., complete topological rings $S$ containing an open subring $S_0 \subseteq S$ on which the topology is $\pi$-adic for some $\pi \in S_0$ such that $S = S_0[\frac{1}{\pi}]$. Such an element $\pi$ is called \emph{pseudo-uniformizer}. A subset $M$ of $S$ is then called \emph{bounded}, if there exist an $n \geq 1$ such that $\pi^nM \subseteq S_0$, a condition that is independent of the choice of $S_0$ and $\pi$. We denote by $S^{\circ}$ the subset of powerbounded elements, i.e. of elements $a \in S$ such that the set $\set{a^n}{n \in \NN}$ is bounded. Then $S^{\circ}$ is a subring of $S$.

A complete Tate ring $S$ is called \emph{perfectoid} if $S^{\circ}$ is bounded and $S^{\circ}$ is perfectoid in the above sense. This definition is shown to be equivalent to Fontaine's definition in \cite[3.20]{BMS_IntegralHodge}. To define the tilt for a perfectoid complete Tate ring $S$, one chooses a pseudo-uniformizer $\pi$ such that $\pi^p$ divides $p$ in $S^{\circ}$ and that admits a compatible\footnote{``compatible'' means that $(\pi^{1/p^{n+1}})^p = \pi^{1/p^n}$ for all $n$} sequence of $p$-th power roots $\pi^{1/p^n}$ in $S^{\circ}$ (such a $\pi$ always exists by \cite[6.2.2]{SW_Berkeley}). Modulo $p$ this sequence defines an element $\pi^{\flat} = (\pi,\pi^{1/p},\dots) \in (S^{\circ})^{\flat}$. Then $S^{\flat} := (S^{\circ})^{\flat}[1/\pi^{\flat}]$ is a perfect complete Tate ring, called the \emph{tilt of $S$}. One has $(S^{\flat})^{\circ} = (S^{\circ})^{\flat}$, $\pi^{\flat}$ is a pseudo-uniformizer of $S^{\flat}$, and $S^{\flat\circ}/\pi^{\flat} \cong S^{\circ}/\pi$.

\begin{example}\label{PerfectoidFields}
Every non-archimedean field $(K,|\cdot|)$ is a complete Tate ring. One can take as pseudo-uniformizer any element $\pi$ with $0 < |\pi| < 1$. A subset $M \subseteq K$ is then bounded if and only if $|M|$ is a bounded subset of $\RR_{\geq0}$. Hence the ring of powerbounded elements is the ring of integers $S^{\circ} = O_K = \set{a \in K}{|a| \leq 1}$. Then $K$ is perfectoid if and only if $K$ is not discretely valued (i.e., $|K^{\times}|$ is not a discrete subgroup of $\RR_{>0}$), $|p| < 1$, and the Frobenius $O_K/p \to O_K/p$ is surjective (\cite[3.8]{Scholze_EtCohDiamond}). Any perfectoid complete Tate ring that is a field is of this form by a result of Kedlaya \cite{Kedlaya_CommNonArchBanachFields}. The tilt of a perfectoid field is again a perfectoid field.

\begin{assertionlist}
\item
The non-archimedean field $\QQ_p$ is not perfectoid because there exists no $\pi \in \ZZ_p$ such that $\pi^p$ divides $p$.
\item
Let $C$ be an algebraically closed non-archimedean field extension of $\QQ_p$ (e.g., the completion of an algebraic closure of $\QQ_p$). Then $C$ is perfectoid. Its tilt is an algebraically closed non-archimedean field $C^{\flat}$ of characteristic $p$.
\item
Choose in an algebraic closure of $\QQ_p$ for all $n \geq 1$ a primitive $p^n$-th root of unity $\zeta_n$. Then the completion $\QQ_p^{\rm cycl}$ of $\QQ_p(\{\zeta_n\}_{n\geq 1})$ is a perfectoid field. Its tilt is the $t$-adic completion of the perfect field $\FF_p(\{t^{1/p^n}\}_{n\geq1})$.
\item
If $K$ is any non-archimedean local field with pseudo-uniformizer $\pi$, then the completion of the algebraic extension $K[\pi^{1/p^{\infty}}]= K[\set{\pi^{1/p^n}}{n \geq 1}]$ is a perfectoid field.
\end{assertionlist}
\end{example}

\begin{example}\label{PerfectoidCharp}
A ring (resp.~a complete Tate ring) $S$ of characteristic $p$ is perfectoid if and only if it is perfect (\cite[3.15]{BMS_IntegralHodge} and \cite[3.5]{Scholze_EtCohDiamond}). In this case one has $S = S^{\flat}$.
\end{example}


\subsection{Perfectoid spaces and tilting}

\subsubsection{Adic spaces}

Just as one geometrizes rings by defining their spectrum and then globalizes to schemes in classical Algebraic Geometry, one can do a similar geometrization/globalization process with perfectoid Tate rings. Scholze used for this the notion of adic spaces that was developed in the 1990s by R.~Huber (\cite{Huber_ContVal}, \cite{Huber_FormalRigid}, see also the notes \cite[II--V]{SW_Berkeley}, \cite{Conrad_Adic}, \cite{Morel_Adic}, \cite{Wd_Adic}). They are by definition geometric objects that are locally the adic spectrum of a certain quite general class of topological rings, which include the class of complete Tate rings.

In fact, the adic spectrum $X = \Spa(R,R^+)$ is attached to a pair consisting of a so-called \emph{Huber ring} $R$, which one can always assume to be complete, and a certain open subring $R^+$, called a \emph{ring of integral elements}. There is always a largest ring of integral elements, namely the subring $R^{\circ}$ of powerbounded elements in $R$. One sets $\Spa R := \Spa(R,R^{\circ})$. Although $R^+$ plays an important technical role in the theory, we will sometimes ignore it. 

The adic spectrum carries a presheaf $\Oscr_X$ of complete topological rings which is in fact a sheaf in all cases of interest to us. The global sections of $\Oscr_X$ is the ring $R$. The (pre)sheaf of functions bounded by $1$ is denoted by $\Oscr_X^+$, and its global sections is the ring $R^+$. Adic spaces that are isomorphic to adic spectra are called \emph{affinoid}. A morphism $f\colon X \to Y$ of adic spaces is locally on source and target of the form $\Spa(S,S^+) \to \Spa(R,R^+)$ which is given by a continuous ring homomorphism $\alpha\colon R \to S$ with $\alpha(R^+) \subseteq S^+$.

Adic spaces that are locally the adic spectrum of a complete Tate ring are called \emph{analytic}. The following examples show that the category of adic spaces is quite large.%
\begin{assertionlist}
\item
Let $K$ be a non-archimedean field. Then one can attach to every variety $X$ over $K$ (or, more generally, to every rigid analytic space over $K$) an analytic adic space $X^{\rm ad}$ and the functor $(\cdot)^{\rm ad}$ is fully faithful.
\item
There is a fully faithful embedding from the category of locally noetherian formal schemes to the category of adic spaces.
\item
There is a fully faithful embedding from the category of all schemes to the category of adic spaces.
\end{assertionlist}

\subsubsection{Perfectoid spaces}

We can now give the definition of one of the central notions in Scholze's work.

\begin{definition}\label{DefPerfectoid}
Adic spaces that are locally the adic spectrum of a perfectoid Tate ring are called \emph{perfectoid spaces}.
\end{definition}

The tilting functor glues to a functor $X \sends X^{\flat}$ from the category of perfectoid spaces to the category of perfectoid spaces in characteristic $p$. 

Before coming to the main theorem about perfectoid spaces and tilting one has to agree what it means to be \'etale for a morphism of perfectoid spaces. As perfectoid spaces are always reduced, one cannot use the usual definition via an infinitesimal lifting criterion if one wants to stay within the perfectoid world. But it is easy to define the notion of a finite \'etale morphism. Then one uses that for rigid analytic varieties (or more general for noetherian analytic adic spaces), where one can define ``\'etale'' by unique infinitesimal lifting, a morphism is \'etale if and only if it can be factorized locally on target and source into an open immersion followed by a finite \'etale morphism\footnote{The analogue assertion does not hold for schemes (not even for \'etale morphisms of curves over a field).}. Hence one ends up with the following definition.

\begin{definition}(\cite[Section 7]{Scholze_Perfectoid})\label{DefEtale}
Let $f\colon X \to Y$ be a morphism of perfectoid spaces.%
\begin{assertionlist}
\item
The morphism $f$ is called \emph{finite \'etale} if for every open affinoid perfectoid subspace $V = \Spa(R,R^+) \subseteq Y$, the preimage $U = f^{-1}(V) = \Spa(S,S^+) \subseteq X$ is affinoid perfectoid, $R \to S$ is a finite \'etale morphism of rings, and $S^+$ is the integral closure of $R^+$ in $S$.
\item
The morphism $f$ is called \emph{\'etale} if it can be locally on $X$ and $Y$ be factorized into an open immersion followed by a finite \'etale morphism.
\end{assertionlist}
\end{definition}

Now the fundamental theorem about perfectoid spaces is as follows.

\begin{theorem}(\cite[3.18,3.20, 3.24]{Scholze_EtCohDiamond})\label{ThmPerfd}
Let $X$ be a perfectoid space and let $X^{\flat}$ be its tilt.
\begin{assertionlist}
\item\label{ThmPerfd1}
The tilting functor $Y \sends Y^{\flat}$ from perfectoid spaces over $X$ to perfectoid spaces over $X^{\flat}$ is an equivalence of categories.
\item\label{ThmPerfd2}
Every adic space that is \'etale\footnote{With the definition of ``\'etale'' as in Definition~\ref{DefEtale}} over $X$ is perfectoid and $Y \sends Y^{\flat}$ induces an equivalence between the \'etale site\footnote{One defines the notion of the \'etale site of a perfectoid space using Definition~\ref{DefEtale} as one defines the \'etale site as for schemes, see the introduction.} of $X$ and the \'etale site of $X^{\flat}$. In particular the categories of \'etale sheaves on $X$ and on $X^{\flat}$ are equivalent.
\item
If $X$ is the adic spectrum of a perfectoid Tate ring $R$, then there is the following acyclicity result: One has $H^i(X,\Oscr_X) = 0$ for all $i > 0$ and $\pi H^i(X,\Oscr_X^+) = 0$ for all $i > 0$ and for every pseudo-uniformizer $\pi$ of $R$.
\end{assertionlist}
\end{theorem}

Part~\ref{ThmPerfd2} of the theorem (and its proof) gives simultaneously a vast generalization of the result of Fontaine and Wintenberger and a strong form of Faltings' almost purity theorem, which was one of Faltings' key techniques in his work on p-adic Hodge theory.

Note that without fixing a perfectoid base $X$ and its tilt $X^{\flat}$ the functor $Y \sends Y^{\flat}$ is not ``injective''. In fact the ``moduli space'' of all untilts of a perfectoid space $S$ in characteristic $p$ is the Fargues-Fontaine curve attached to $S$ (see Subsection~\ref{FFCURVE} below).

Let us conclude with an example about ``naive perfectoidization'', i.e., attaching to a ``classical'' variety a perfectoid space. Say we are given a variety $X$ over a non-archimedean local field $K$ of characteristic $0$. Often one can check properties of $X$ after passing to a field extension which allows to replace $K$ by a perfectoid field extension. But even then $X$ will not be perfectoid (except if $X$ is zero-dimensional). There is however the following construction of perfectoidization in a very special case.

\begin{example}\label{PerfectoidToric}
Let $K$ be a perfectoid field. Let $\sigma \subseteq \RR^n$ be a finitely generated rational polyhedral cone and let $K[\sigma \cap \ZZ^n]$ be the monoid algebra over $K$ attached to the abelian monoid $\sigma \cap \ZZ^n$. It is a finitely generated $K$-algebra. We consider the $K$-scheme $X := \Spec K[\sigma \cap \ZZ^n]$. For instance, if $\sigma \cap \ZZ^n = \NN_0^n$, then $X = \AA^n_K$ is the $n$-dimensional affine space of $K$.

Now $K\langle \sigma \cap \ZZ^n \rangle$ is a complete Tate ring. It consists of power series whose coefficients are indexed by $\sigma \cap \ZZ^n$ and that converge to $0$ if their index goes to infinity. We denote the attached adic spectrum by $\Xcal$\footnote{As noted above, one also has to specify a ring of integral elements. But in this case there is a canonical choice, namely $O_K\langle \sigma \cap \ZZ^n \rangle$. A similar remark holds for the perfectoidization.}. Note that this is often not the adic space $X^{\rm ad}$ attached to $X$. For instance it attaches to the affine space the $n$-dimensional closed unit ball of dimension $n$. We can now attach to this adic space a perfectoid space. It is defined as the adic spectrum of $K\langle \sigma \cap \ZZ^n[1/p] \rangle$.

This construction can be generalized to arbitrary toric varieties because these are glued from schemes of the form $\Spec K[\sigma \cap \ZZ^n]$. For instance, one obtains a perfectoidization of the projective space. Moreover, suppose that $X^{\rm perfd}_K$ and $X^{\rm perfd}_{K^{\flat}}$ are the constructed perfectoidizations of toric varieties attached to the same collection of cones, one over $K$ and the other one over its tilt $K^{\flat}$. Then the tilt of $X^{\rm perfd}_K$ is isomorphic to $X^{\rm perfd}_{K^{\flat}}$ by \cite[8.5]{Scholze_Perfectoid}.
\end{example}

There is also a less naive way of perfectoidization (see \cite[8]{BS_Prismatic}) that functorially attaches to every $p$-adically complete algebra $S$ over a perfectoid ring $R$ its perfectoidization $S_{\perfd}$. If $S$ is of characteristic $p$ it is the colimit perfectization $\colim_{x\sends x^p}S$. In general, $S_{\perfd}$ is not a perfectoid ring in the above sense but an object of higher algebra, more precisely a $p$-complete $\EE_{\infty}$-$S$-algebra (see Subsection~\ref{PERFDIZATION} for some details).


\section{First applications of perfectoid spaces}

There are a number of applications of the direct passage from mixed (or zero) characteristic to characteristic $p$ given by the tilting process. Each time, the idea is to reduce a question in characteristic zero or mixed characteristic by tilting to a question in positive characteristic, where the question was already answered. We give two examples.

\subsection{The weight monodromy conjecture}\label{WEIGHTMONODROMY}

The first one, given by Scholze in his thesis \cite{Scholze_Perfectoid}, was a proof of Deligne's weight monodromy conjecture for a large class of algebraic varieties. Let us briefly recall this conjecture.

Let $K$ be a non-archimedean local field (i.e., $K$ is a finite extension of $\QQ_p$ or of $\FF_p((t))$). Let $\Kbar$ be an algebraic closure of $K$. Let $O_K$ be its ring of integers and let $O_{\Kbar}$ be the integral closure of $O_K$ in $\Kbar$. Then $O_K$ and $O_{\Kbar}$ are local rings, the residue field of $O_K$ is a finite field $k$ of characteristic $p$ and the residue field of $O_{\Kbar}$ is an algebraic closure $\kbar$ of $k$. Every element in the absolute Galois group $\Gamma_K := \Gal(\Kbar/K)$ induces by restriction to $O_{\Kbar}$ an automorphism of $O_{\Kbar}$ fixing $O_K$. Modulo the maximal ideal of $O_{\Kbar}$ one obtains an automorphism of $\kbar$ fixing $k$, i.e., an element of $\Gamma_k := \Gal(\kbar/k)$. This defines a surjective group homomorphism $\Gamma_K \to \Gamma_k$. In $\Gamma_k$ there is the so-called \emph{geometric Frobenius} given by $x \sends x^{1/q}$, where $q$ is the number of elements of $k$. It is the inverse of the (arithmetic) $q$-Frobenius $x \sends x^q$. Any element $\phi$ in $\Gamma_K$ that maps to the geometric Frobenius in $\Gamma_k$ is again called a geometric Frobenius.

Now let $X$ be a proper smooth variety over $K$ and let $H^i(X_{\Kbar},\QQbar_{\ell})$ be its $\ell$-adic cohomology (see also Subsection~\ref{PROETSCHEME} below). It is a finite-dimensional vector space over an algebraic closure of $\QQ_{\ell}$, where $\ell$ is a prime number different from $p$. It should be considered as the algebraic-geometric variant of (rational) singular cohomology for manifolds. It is endowed with an action of $\Gamma_K$, and the structure of these $\Gamma_K$-modules $H^i(X_{\Kbar},\QQbar_{\ell})$ is a central invariant of $X$.

If $X$ has good reduction, i.e., there exists a smooth proper model $\Xscr$ over $O_K$ (see Subsection~\ref{NONARCH}), then the action of $\Gamma_K$ on $H^i(X_{\Kbar},\QQbar_{\ell})$ factors through $\Gamma_k$ and one has an isomorphism of $\Gamma_k$-modules $H^i(X_{\Kbar},\QQbar_{\ell}) \cong H^i(\Xscr_{\kbar},\QQbar_{\ell})$, where the right hand side is the $\ell$-adic cohomology of the special fiber of $\Xscr$. As $\Gamma_k$ is generated as a profinite group by the geometric Frobenius, the $\Gamma_k$-action is given by the single endomorphism induced by the action of the geometric Frobenius $F_q$. Its characteristic polynomial is determined by counting the points $\Xscr_k(k_n)$, where $k_n$ is the unique extension of $k$ of degree $n$ within $\kbar$ by \eqref{EqGrothLefschetz}\footnote{Indeed, a standard argument (e.g., \cite[Rapport 3.1]{SGA4}) shows that $T\frac{d}{dT}\log(\det(1 -F_q T)^{-1}) = \sum_{n \geq 1}\tr(F_q^n | H^i(\Xscr_{\kbar},\QQbar_{\ell}))$.}. Moreover, by Deligne's proof of the Weil conjectures (\cite{Deligne_Weil1}, \cite{Deligne_Weil2}), the characteristic polynomial $\det(T - F_q | H^i(\Xscr_{\kbar},\QQbar_{\ell}))$ has coefficients in $\ZZ$ that are independent of $\ell$. All complex roots $\alpha$ of this polynomial (i.e., the complex eigenvalues of $F_q$ on $H^i(\Xscr_{\kbar},\QQbar_{\ell})$) have absolute value $q^{i/2}$. One says, that $F_q$ has weight $i$ on $H^i(\Xscr_{\kbar},\QQbar_{\ell})$.

In general, $X$ not necessarily has good reduction. Then one important piece still missing for $H^i(X_{\Kbar},\QQbar_{\ell})$ to have a similar description is Deligne's weight monodromy conjecture. By Grothendieck's monodromy theorem, there exists a nilpotent operator $N$, called the \emph{monodromy operator}, on $V = H^i(X_{\Kbar},\QQbar_{\ell})$\footnote{It is characterized by the property that the restriction of the $\Gamma_K$-action to an open subgroup (automatically of finite index) of $I := \Ker(\Gamma_K \to \Gamma_k)$ is given by $\exp(Nt_{\ell})$, where $t_{\ell}\colon I \to \ZZ_{\ell}(1)$ is the maximal pro-$\ell$-quotient, given by the limit of the homomorphism $t_{\ell,n}\colon I \to \mu_{\ell^n}$ defined by choosing a system of $\ell^n$-th roots $\varpi^{1/\ell^n}$ of a uniformizing element $\varpi$ of $K$ and requiring $\sigma(\varpi^{1/\ell^n}) = t_{\ell,n}(\sigma)\varpi^{1/\ell^n}$ for $\sigma \in I$.}. If $X$ has good reduction, then $N = 0$. The action of $\Gamma_K$ on $V$ is characterized by the action of a geometric Frobenius and the monodromy operator. An elementary linear algebra argument shows the existence of a unique increasing filtration, called the \emph{monodromy filtration}, by subspaces $\Fil_j(V)$ such that $N(\Fil_j(V)) \subseteq \Fil_{j-2}(V)$ and $N^j$ induces an isomorphism $\gr_jV := \Fil_j(V)/\Fil_{j-1}(V) \to \gr_{-j}(V)$ for all $j \geq 0$. Now Deligne's conjecture predicts that for any geometric Frobenius $\phi$ in $\Gamma_k$ all eigenvalues $\alpha$ of $\phi$ in $\gr_jV$ are algebraic and $|\alpha| = q^{(i+j)/2}$ for all complex absolute values.

If $K$ is a local field of characteristic $p$, then the conjecture is essentially known\footnote{``Essentially known'' means that the conjecture is known for varieties that are obtained from varieties over a global field of characteristic $p$ whose completion at some place is $K$. This is a result by Deligne himself.}. For local fields of characteristic $0$ this was previously known in dimension $\leq 2$ by work of Rapoport, Zink, and de Jong. Scholze proved as an application of his theory of perfectoid spaces the following theorem in his thesis.

\begin{theorem}\label{Monodromy}
Let $K$ be a non-archimedean local field of characteristic $0$ and let $X$ be a geometrically connected proper smooth scheme that is a set-theoretic complete intersection in a projective smooth toric variety. Then the weight monodromy conjecture holds for $X$.
\end{theorem}

The theorem is reduced to Deligne's theorem via tilting. One first replaces $K$ by a perfectoid extension. If $X$ is a smooth proper toric variety, then the theorem can easily be deduced from Example~\ref{PerfectoidToric}. In general, if $X$ is a closed subvariety of a smooth projective toric variety $Z$, one needs to approximate the pre-image of $X^{\rm ad}$ under the homeomorphism $Z^{\perfd}_{K^{\flat}} \to Z^{\perfd}_K$ by an algebraic object. This is non-trivial and here the hypothesis that $X$ is set-theoretically a complete intersection in $Z$ and that the ambient toric variety is projective is used.


\subsection{Hochster's conjecture}

One of the most important conjectures in commutative algebra was Hochster's conjecture that deals with the question whether given a noetherian ring $R$ and a finite ring extension $R \subseteq S$, $R$ is always a direct summand of $S$ as an $R$-module. Hochster conjectured around 1969 that this is the case if $R$ is regular\footnote{In general, $R$ is not a direct summand as the example $R = \QQ[X,Y]/(XY)$ and $S$ its normalization shows: the map $R \to S$ does not stay injective modulo $X + Y$.}. This conjecture is (non-trivially) equivalent to a number of other important conjectures in commutative algebra. Examples are the following.
\begin{assertionlist}
\item
Given an injective integral map $A \to B$ of noetherian rings and an $A$-module $M$. If $M \otimes_A B$ is a flat $B$-module, then $M$ is a flat $A$-module.
\item
An $r$-th module of syzygies of a finitely generated module over a regular local ring, if not free, has torsion-free rank at least $r$.
\item
Let $R$ be a complete local noetherian domain of Krull dimension $d$, let $\mfr$ be its maximal ideal, and let $\Rbar$ be the integral closure of $R$ in an algebraic closure of the field of fractions of $R$. Then $H^d_{\mfr}(\Rbar) \ne 0$.
\end{assertionlist}

Hochster himself proved his conjecture for rings containing a field\footnote{If this field has characteristic $0$, then this is an easy trace map argument that shows that the conjecture holds in this case more generally if $R$ is only assumed to be normal.} and showed that to prove the conjecture in general it is sufficient to treat the case where $R$ is a complete local noetherian ring of mixed characteristic $(p,0)$. This was settled by Andr\'e in \cite{Andre_Hochster} by reducing the case of mixed characteristic to the case of characteristic $p$ (in which case that ring contains a field) by using perfectoid rings and tilting. Hence all of the above conjectures are now theorems.


\section{The pro-\'etale topology}\label{PROET}

The \'etale topology for schemes was invented by Grothendieck to construct a cohomology theory for varieties over arbitrary fields to prove the Weil conjectures. But it has many further applications within and outside Arithmetic Geometry, for instance the construction of representations of finite groups of Lie type by Deligne and Lusztig. The \'etale topology is a Grothendieck topology (see the introduction) and can be considered as a refinement of the Zariski topology. Moreover, over the complex numbers it can be viewed (with a grain of salt) as an algebraic variant of the complex analytic topology. For instance, a morphism of complex varieties is \'etale if and only if it is a local isomorphism if considered as a holomorphic map of complex analytic spaces.

Scholze extends this theory with strong consequences. The main focus lies on adic and perfectoid spaces but even for usual schemes one gets new results. Let us consider this case first.


\subsection{The pro-\'etale topology for schemes}\label{PROETSCHEME}

Bhatt and Scholze defined and studied in \cite{BS_ProEt} the pro-\'etale topology for schemes. It is the Grothendieck topology ``generated'' by the Zariski topology and surjective morphisms of the form $\Spec A \to \Spec R$, where $A$ is a filtered colimit of \'etale $R$-algebras\footnote{One then says that $\Spec A \to \Spec R$ is a pro-\'etale morphism of schemes. But this notion of pro-\'etale is badly behaved (for instance it is not local on the target) and hence it is often better to work with the notion of being weakly \'etale (\cite[Tag 092A]{Stacks}). Bhatt and Scholze show in \cite{BS_ProEt} that both notions yields the same theory of sheaves and the same cohomology groups.}. It is therefore finer than the \'etale topology and coarser than the fpqc-topology. Roughly spoken, pro-\'etale cohomology has all the advantages of the \'etale topology but not (some of) its problems.

\'Etale cohomology on a scheme $X$ works best for abelian torsion sheaves that are killed by some integer $n \geq 1$ that is invertible on $X$. To define a cohomology theory with coefficients in a field of characteristic zero (as needed for instance to have a good Lefschetz fixed-point formula as in Subsection~\ref{COUNTPOINT}) one chooses a prime $\ell$ that is invertible on $X$ and then the classical definition of $\ell$-adic cohomology is
\[
H^i_{\rm et}(X,\QQ_{\ell}) := (\lim_n H^i_{\rm et}(X,\ZZ/\ell^n\ZZ)) \otimes_{\ZZ_{\ell}} \QQ_{\ell}
\]
to obtain a vector space over $\QQ_{\ell}$. Note that one cannot (and does not want to) replace $\lim_n H^i_{\rm et}(X,\ZZ/\ell^n\ZZ)$ by $H^i_{\rm et}(X,\ZZ_{\ell})$. Although the latter makes formally sense, it has not the desired properties (for instance it is always zero for $i = 1$ and $X$ normal). But with this definition, the $\ell$-adic cohomology does not arise as derived functors. This causes problems with functoriality, for instance for the long exact cohomology sequence associated to a short exact sequence of $\ell$-adic sheaves. Some of these problems have been overcome by modifying the definition, notably by Jannsen and Ekedahl.

Bhatt and Scholze had the insight that a solution (encompassing the settings of Jannsen and Ekedahl) would be to replace the \'etale topology by the pro-\'etale topology. For the pro-\'etale topology the usual cohomology with values in the sheaves attached to the topological rings $\ZZ_{\ell}$ and $\QQ_{\ell}$ has all the desired properties if $\ell$ is invertible on the scheme\footnote{It remains the problem that cohomology with values in $p$-torsion sheaf or in $\ZZ_p$, where $p$ is a prime not invertible on the scheme. This is of principal nature, see the section on prismatic cohomology below.}. The pro-\'etale $\ell$-adic cohomology thus is a derived functor and satisfies all good properties of those. They also showed that for torsion sheaves annihilated by some integer invertible on the scheme (i.e., for those sheaves for which the \'etale topology works just fine), \'etale and pro-\'etale cohomology coincide. Therefore, to put it somewhat polemically, one can forget nowadays about \'etale cohomology and work instead with pro-\'etale cohomology only.

There are aspects that give the pro-\'etale topology a very different flavor than the \'etale topology. For instance, the topology is sufficiently fine such that for every quasi-compact scheme $X$ there exists a covering $Z \to X$ for the pro-\'etale topology, where $Z$ is a ``very disconnected'' affine scheme. More precisely (\cite[2]{BS_ProEt}), every connected component of $Z$ has a unique closed point, the subspace $Z^c$ of closed points is closed in $Z$ and extremally disconnected\footnote{A topological space $S$ is called \emph{extremally connected} if it is compact Hausdorff and the closure of every open set is open. Equivalently, $S$ is a projective object in the category of compact Hausdorff spaces.}. Moreover, one can assume that for every closed point $z \in Z$ the local ring $\Oscr_{Z,z}$ is strictly henselian.

Such covering spaces $Z$ are almost never noetherian, even if $X$ was noetherian\footnote{In fact, such an affine scheme $Z = \Spec A$ is noetherian if and only if $A$ is the product of a finite number of strictly henselian noetherian local rings. Then $Z^c$ is a finite discrete space}. Hence one has coverings by objects that are ``almost discrete'' but ``very non-finite''. This is a typical phenomenon which will also hold for the pro-\'etale topology of perfectoid spaces which is described next.


\subsection{The pro-\'etale topology and the v-topology for perfectoid spaces}\label{PROETPERFD}

A morphism of perfectoid spaces is called \emph{pro-\'etale} if locally on source and target it is a morphism $f\colon X \to Y$ of affinoid perfectoid spaces, where $X$ can be written as a cofiltered inverse limit\footnote{One has to be careful about the precise meaning of the inverse limit. It is {\em not} the inverse limit in the category of adic spaces. See \cite[6.4]{Scholze_EtCohDiamond} or \cite[4.1]{SW_ModuliPDiv} for a precise definition.} of \'etale maps of affinoid perfectoid spaces $X_i \to Y$. As usual this yields the pro-\'etale topology on the category $\Perfd$ of perfectoid spaces. Scholze also defines a much finer topology on $\Perfd$, the $v$-topology where essentially all surjective families of morphisms of perfectoid spaces are defined to be coverings\footnote{All coverings are always defined to be subject to a standard quasi-compactness assumption that ensures that one can refine coverings of quasi-compact spaces always by a finite covering of quasi-compact spaces. This condition also appears in the definition of the fpqc topology of schemes. For topologies whose coverings are defined to be surjective families of certain open morphisms (such as the \'etale or the fppf topology) this condition is automatic.}, see \cite[8.1]{Scholze_EtCohDiamond} for the precise definitions.

Clearly every sheaf for the v-topology is in particular a sheaf for the pro-\'etale topology. It is somewhat amazing that such a fine topology as the v-topology is useful. But in fact Scholze shows:

\begin{theorem}(\cite[8.7, 8.8]{Scholze_EtCohDiamond})\label{ProEtPerfectoid}
\begin{assertionlist}
\item
Presheaves represented by perfectoid spaces are sheaves for the $v$-topology. The structure presheaves $\Oscr\colon X \sends \Oscr_X(X)$ and $\Oscr^+\colon X \sends \Oscr^+_X(X)$ are sheaves for the v-topology.
\item
Let $X = \Spa(R,R^+)$ be an affinoid perfectoid space. Then for all $i > 0$ one has $H^i_v(X,\Oscr) = H^i(X_{\proet},\Oscr) = 0$ and $\varpi H^i_v(X,\Oscr^+) = \varpi H^i(X_{\proet},\Oscr^+) = 0$ for every pseudo-uniformizer $\varpi \in R^+$.
\end{assertionlist}
\end{theorem}

Scholze shows this amazing result by first proving all assertions for the pro-\'etale topology. Then he shows that every perfectoid space has a pro-\'etale universally open covering by affinoid perfectoid spaces $Z$ that are \emph{strictly totally disconnected}, i.e., every connected component of $Z$ is of the form $\Spa(C,C^+)$ where $C$ is an algebraically closed perfectoid field and $C^+ \subseteq C$ is an open and bounded valuation ring (\cite[7.18]{Scholze_EtCohDiamond}). Adic spaces of the form $\Spa(C,C^+)$ are of a particular simple form. They have a unique closed point and their topological spaces are totally ordered by specialization of points. Then to obtain results for the v-topology, Scholze shows that every surjective morphism from a perfectoid space to a strictly totally disconnected spaces is automatically faithfully flat (\cite[7.23]{Scholze_EtCohDiamond}).

\subsection{The pro-\'etale site of an analytic adic space}

The pro-\'etale site $X_{\textup{pro\'et}}$ can also be defined for general analytic adic spaces $X$, for instance if $X$ is the adic space attached to a rigid analytic variety over a local field of characteristic $0$. Then by \cite[15.3]{Scholze_EtCohDiamond} (see also \cite[8.1]{SW_Berkeley}) one has:

\begin{proposition}\label{LocallyPerfectoid}
Let $X$ be an analytic adic space over $\ZZ_p$. Then $X$ has a pro-\'etale covering $\Xtilde \to X$ with $\Xtilde$ perfectoid.
\end{proposition}

For instance, if $X = \Spa \QQ_p$, then one can choose $\Xtilde = \Spa \QQ_p^{\rm cycl}$ (Example~\ref{PerfectoidFields}).


\section{Scholze's way to pass to positive characteristic II: Diamonds}\label{DIAMOND}

\subsection{Definition of diamonds}

The basic references are \cite{SW_Berkeley}, which provides also a lot of ideas and examples, and \cite{Scholze_EtCohDiamond}.

Spaces occurring ``in nature'' such as varieties over non-archimedean local fields are rarely perfectoid. If one wants to extend the tilting functor to such spaces to obtain a functor $(\cdot)^{\diamond}$ one could proceed roughly as follows.

Using Proposition~\ref{LocallyPerfectoid} write your space $X$ that you are interested in as $\Xtilde/R$, where $\Xtilde$ is perfectoid and $R = \Xtilde \times_X \Xtilde$ is a perfectoid equivalence relation and then set $X^{\diamond} := \Xtilde^{\flat}/R^{\flat}$. Of course, there are number of problems with this approach. For instance, one has to find a category in which $\Xtilde^{\flat}/R^{\flat}$ makes sense. And then one has to check that $X^{\diamond}$ is independent of the choice if $\Xtilde$ and $R$.

Another approach is to note that if $X$ is a perfectoid space, then $X^{\flat}$ represents by Theorem~\ref{ThmPerfd}~\ref{ThmPerfd1} the functor $X^{\diamond}$ on perfectoid spaces $S$ of characteristic $p$ that sends $S$ to the set of isomorphism classes of \emph{untilts}, i.e., perfectoid spaces $S^{\sharp}$ with an isomorphism $\iota\colon (S^{\sharp})^{\flat} \cong S$, together with a morphism $S^{\sharp} \to X$ of perfectoid spaces. This functor then also makes sense for more general $X$.

Scholze defines in \cite{Scholze_EtCohDiamond} (see also \cite{SW_Berkeley}) a category in which $X^{\diamond}$ lives, the category of diamonds, and shows that both approaches are valid and give the same object. The second approach suggests for the category of diamonds to be a category of certain sheaves on the category $\Perf$ of all perfectoid spaces in characteristic $p$ (considered as a full subcategory of $\Perfd$). The first approach suggests that these sheaves should be quotients of a perfectoid space by some pro-\'etale equivalence relation. Hence we arrive at the following definition.

\begin{definition}\label{DefDiamond}
A \emph{diamond} is a sheaf $Y$ for the pro-\'etale topology on $\Perf$ such that $Y$ can be written as a quotient $X/R$, where $X$ is representable by a perfectoid space and $R$ is a pro-\'etale equivalence relation on $X$, i.e., a representable equivalence relation $R \subseteq X \times X$ such that the projections $s,t\colon R \to X$ are pro-\'etale.
\end{definition}

The definition is reminiscent of that of an algebraic space which is defined, roughly spoken, as the quotient of a scheme by an \'etale equivalence relation. But as the pro-\'etale topology is much finer than the \'etale topology, this analogy has to be taken with a grain of salt.

Scholze shows that every diamond is automatically also a sheaf for the $v$-topology (\cite[11.9]{Scholze_EtCohDiamond}). This is an analogy to a result of Gabber that shows that every algebraic space is a sheaf for the fpqc topology (\cite[0APL]{Stacks}). The category of diamonds has products (\cite[8.3.7]{SW_Berkeley}) and fiber products (\cite[11.4]{Scholze_EtCohDiamond}). If $Y = X/R$, one defines its underlying topological space by $|Y| := |X|/|R|$ (endowed with the quotient topology). The topological space $|Y|$ is independent of the presentation $Y = X/R$ (\cite[11.13]{Scholze_EtCohDiamond}). One defines the notion of an \'etale (resp.~finite \'etale) map of diamonds (\cite[10.1]{Scholze_EtCohDiamond}).


\subsection{Extending the tilting functor}

As suggested above one can now attach to every analytic space $X$ over $\ZZ_p$ its diamond $X^{\diamond}$. It is the v-sheaf on $\Perf$ defined by
\[
X^{\diamond}(S) = \{(S^{\sharp},\iota,f\colon S^{\sharp} \to X)\}/\cong
\]
where $(S^{\sharp},\iota)$ is an untilt of $S$ and $f$ is a morphism of adic spaces. If $X$ is perfectoid, then $X^{\diamond}$ is represented by $X^{\flat}$. If $X = \Spa(A,A^+)$ is affinoid, we also write $\Spd(A,A^+)$ instead of $X^{\diamond}$. 

\begin{theorem}(\cite[15.6]{Scholze_EtCohDiamond}, \cite[10.2.3]{SW_Berkeley})\label{PropertyDiamand}
Let $X$ be an analytic space over $\ZZ_p$. Then $X^{\diamond}$ is a diamond. Moreover:
\begin{assertionlist}
\item
The underlying topological spaces $|X|$ and $|X^{\diamond}|$ are homeomorphic. One has equi\-valences of \'etale sites $X^{\diamond}_{\et} \cong X_{\et}$ and finite \'etale sites $X^{\diamond}_{\textup{f\'et}} \cong X_{\textup{f\'et}}$.
\item
Let $K$ be a non-archimedean field over $\QQ_p$. Then the functor $X \sends X^{\diamond}$ from the category of normal\footnote{In fact, the diamond functor is also fully faithful on the category of seminormal rigid analytic spaces as defined in \cite[3.7]{KedlayaLiu_II}. This is useful since the pro-\'etale site of a rigid analytic space and of its semi-normalization are isomorphic by \cite[8.2.2]{KedlayaLiu_II}.} rigid analytic spaces over $K$ to the category of diamonds over $\Spd(K,O_K)$ is fully faithful.
\end{assertionlist}
\end{theorem}


\subsection{\'Etale cohomology of diamonds}

Using these concepts, Scholze has established an \'etale cohomology theory of diamonds in \cite[16ff]{Scholze_EtCohDiamond}. In particular, he constructs the analogue of the ``six-functor calculus'' and appropriate versions of the proper and smooth base change theorems. 

This theory is one of the key tools in the geometric construction of a local Langlands correspondence, see Section~\ref{LOCALLANGLANDS} below. 


\part{$p$-adic Hodge theory}

The starting point of classical Hodge theory is the existence of functorial isomorphisms
\[
H^n(X,\ZZ) \otimes_{\ZZ} \CC = H^n(X,\CC) = H^n_{\dR}(X/\CC) = \bigoplus_{i+j=n}H^j(X,\Omega^{i}_X)
\]
for $X$ a compact K\"ahler manifold. Here $H^n(X,\ZZ)$ and $H^n(X,\CC)$ denotes singular cohomology of $X$ with values in $\ZZ$ or $\CC$. Equivalently, it denotes sheaf cohomology of the constant sheaf $\ZZ$ or $\CC$ on $X$. Moreover, $H^n_{\dR}(X/\CC)$ denotes de Rham cohomology of $X$. The functorial isomorphism $H^n(X,\CC) = H^n_{\dR}(X/\CC)$ is called the de Rham isomorphism\footnote{If one views $H^n(X,\CC)$ as cohomology in the constant sheaf, it is simply given by the map of complexes of sheaves of $\CC$-vector spaces $\CC \to \Omega^{\bullet}_X$ which is a quasi-isomorphism by the Poincar\'e lemma for complex manifolds.}. The last equality is the Hodge decomposition theorem. These results also hold in an algebraic version for a proper smooth variety $X$ over $\CC$. In this case, $H^n(X,\ZZ)$ and $H^n(X,\CC)$ denote the singular cohomology (or, equivalently, the sheaf cohomology of the constant sheaves) of the compact smooth manifold attached to $X$.

The existence of the Hodge decomposition follows from the degeneracy of the first hypercohomology spectral sequence of the de Rham complex
\[
E^{ij}_1 = H^j(X,\Omega^i_X) \implies H^{i+j}_{\dR}(X/\CC).
\]
This spectral sequence is called the Hodge-de Rham spectral sequence (or sometimes also the Fr\"olicher spectral sequence). It yields a Hodge filtration on $H^n_{\dR}(X/\CC)$ whose graded spaces are $H^j(X,\Omega^i_X)$ for $i+j = n$. Moreover, with respect to complex conjugation on $H^n(X,\CC) = H^n(X,\RR) \otimes_{\RR} \CC$ one has $\overline{H^j(X,\Omega^{i}_X)} = H^i(X,\Omega^{j}_X)$ which yields a splitting of the filtration and hence the Hodge decomposition.

For algebraic varieties over general fields $k$ there is no attached complex manifold and one replaces the singular cohomology on manifolds with the (pro-)\'etale cohomology of the variety with coefficients in $\ZZ_{\ell}$ or $\QQ_{\ell}$, where $\ell$ is some prime number that is invertible in $k$. The pro-\'etale cohomology also has the advantage of carrying a Galois action if the varieties are defined over non-algebraically closed field. For instance, if a complex smooth proper variety $X$ is already defined over $\QQ$, we obtain a continuous action by the absolute Galois group $\Gamma_{\QQ}$ of $\QQ$ on the $\ZZ_{\ell}$-module $H^i(X_{\proet},\ZZ_{\ell})$, and thus on the $\QQ_{\ell}$-vector space $H^i(X_{\proet},\QQ_{\ell})$.

An important example of such a Galois module is the Tate module. If $K$ is a field, choose an algebraic closure $\Kbar$. Then $\Gamma_K = \Gal(\Kbar/K)$ acts on $\Kbar$. The induced action on the group $\mu_{m} := \set{a \in \Kbar^{\times}}{a^{m} = 1}$ of $m$-th roots of unity makes $\mu_m$ into a $\Gamma_K$-module over $\ZZ/m\ZZ$. If $m$ is invertible in $K$, then the underlying group of $\mu_m$ is isomorphic to $\ZZ/m\ZZ$ by choosing a primitive $m$-th root of unity. Hence for $\ell$ a prime invertible in $K$ one has a $\Gamma_K$-module over $\ZZ_{\ell}$
\[
\ZZ_{\ell}(1) := \lim_n\mu_{\ell^n},
\]
whose underlying $\ZZ_{\ell}$-module is free of rank $1$. As usual this allows to define $\Gamma_k$-modules $\ZZ_{\ell}(j)$ for any $j \in \ZZ$ by setting $\ZZ_{\ell}(j) := \ZZ_{\ell}(1) \otimes_{\ZZ_{\ell}} \dots \otimes_{\ZZ_{\ell}} \ZZ_{\ell}(1)$ ($j$ factors) for $j \geq 0$ and by defining $\ZZ_{\ell}(j)$ to be the dual $\ZZ_{\ell}(-j)\vdual$ for $j < 0$. For any $\Gamma_K$-module $M$ over $\ZZ_{\ell}$ we set $M(j) := M \otimes \ZZ_{\ell}(j)$ and call it the \emph{$j$-th Tate twist} of $M$.

For instance, if $X = \GG_{m,K} = \Spec K[T,T^{-1}]$, then $H^1(X_{\proet},\ZZ/\ell^n\ZZ) = \mu\vdual_{\ell^n}$ and hence $H^1(X_{\proet},\ZZ_{\ell}) = \ZZ_{\ell}(-1)$.

Hodge theory over $p$-adic fields was initiated by Tate. Since then results by Fontaine, Messing, Faltings, Kato, Tsuji, Niziol, Beilinson and others gave a rather complete picture for smooth proper varieties and the relation of their de Rham cohomology and their \'etale cohomology with coefficients over $\QQ_p$ (as $p$-adic field have characteristic $0$, the prime $p$ is invertible in them).

\section{Rational $p$-adic Hodge theory for rigid analytic spaces}\label{RATIONALHODGE}


The techniques developed by Scholze allow to extend some of the results to rigid analytic varieties. Let $C$ be a complete algebraically closed extension of $\QQ_p$.

\begin{theorem}(\cite[13.3]{BMS_IntegralHodge})\label{DegHodgeSS}
Let $X$ be a proper smooth rigid-analytic space over $C$. Then the Hodge-de Rham spectral sequence
\[
E^{ij}_1 = H^j(X,\Omega^i_{X/C}) \implies H^{i+j}_{\dR}(X/C)
\]
degenerates at $E_1$.
\end{theorem}

It is remarkable that, contrary to the complex case, the theorem holds without any K\"ahler type hypothesis on $X$.

The theorem shows that $H^n_{\dR}(X/C)$ has a canonical filtration (the ``Hodge filtration'') whose graded pieces are $H^j(X,\Omega^i_{X/C})$ for $i+j = n$. Over $\CC$, the Hodge filtration is canonical split because its complex conjugate yields a complement. In the $p$-adic case, there does not exist a canonical splitting.

As the $p$-adic analogue of $H^n(X,\CC) = H^n(X,\QQ) \otimes_{\QQ} \CC$ ($X$ complex variety) one may consider for a proper smooth rigid analytic space $X$ the $C$-vector space $H^n(X_{\et},\QQ_p) \otimes_{\QQ_p} C$. In the $p$-adic setting there is no functorial isomorphism between $H^n(X_{\et},\QQ_p) \otimes_{\QQ_p} C$ and $H^n_{\dR}(X/C)$, although the next theorem below in particular implies that they have the same dimension. Scholze obtained from the Leray spectral sequence for the canonical map of sites
\[
\nu\colon X_{\proet} \to X_{\et}
\]
by a non-trivial calculation (see also the sketch of proof below) a second fundamental spectral sequence in $p$-adic Hodge theory which was  new even for elliptic curves:

\begin{theorem}(\cite[3.20]{Scholze_PerfdSurvey}, \cite[13.3]{BMS_IntegralHodge})\label{HodgeTateSS}
Let $X$ be a smooth proper rigid analytic variety over $C$. There is a \emph{Hodge-Tate spectral sequence}
\[
E^{ij}_2 = H^i(X,\Omega^j_X)(-j) \implies H^{i+j}(X_{\et},\ZZ_p) \otimes_{\ZZ_p} C
\]
which degenerates at $E_2$\footnote{Let us remark on the Tate twists. If $X$ is only defined over $C$, then there is no Galois action and hence it is superfluous to include Tate twists. But very often, $C$ is the completion of an algebraic closure $\Kbar$ of a discretely valued field $K$ and $X$ is already defined over $K$. Then $\Gamma_K = \Gal(\Kbar/K)$ acts on $H^{i+j}(X_{\et},\ZZ_p)$. Moreover, the action of $\Gamma_K$ on $\Kbar$ extends by continuity to $C$ making $C$ into $\Gamma_K$-module, and the limit term of the spectral sequence is endowed with the tensor product action. Introducing the Tate twists on the initial terms then makes the whole spectral sequence $\Gamma_K$-equivariant. Moreover, the filtration on the limit term given by the spectral sequence does then split as $\Gamma_K$-modules over $C$ (\cite[1.8]{Scholze_PAdicHodge}).}.
\end{theorem}

To prove the degeneracy of both spectral sequences one first uses a spreading out principle of Conrad and Gabber to realize $X$ as the fiber of a family defined over a discretely valued field (\cite[13.16]{BMS_IntegralHodge}). This allows to reduce to the case that $X$ is defined over a discretely valued field $K$.
Then one uses that every rigid analytic variety is locally for the pro-\'etale topology perfectoid (Proposition~\ref{LocallyPerfectoid}) and relates $H^n(X_{\et},\ZZ_p) \otimes_{\ZZ_p} C$ and the cohomology of the completed structure sheaf $\widehat{\Oscr}$ on $X_{\proet}$. The last step is then to show that $R^j\nu_*\widehat{\Oscr} \cong \Omega^j_X(j)$.

The theorems above are only a small part of results (by Scholze and by many others) in rational $p$-adic Hodge theory that benefited enormously from the theory of perfectoid spaces. Further examples are the extension of the above results to lisse $\ZZ_p$-sheaves instead of $\ZZ_p$ by Scholze in \cite{Scholze_PAdicHodge} or the work of Liu and Zhu on a $p$-adic Riemann-Hilbert correspondence (\cite{LiuZhu_RHCorresp}). 

%
%
%
%
%
%


\section{Integral $p$-adic Hodge theory and prismatic cohomology}

Let $K$ be a finite extension of $\QQ_p$ and let $X$ be a proper smooth scheme over $K$. Let us assume that there exists a good model, i.e. a smooth proper scheme $\Xscr$ over $O_K$ with $\Xscr_K \cong X$. Denote by $k$ the residue field of $O_K$ and by $\Xscr_k$ the special fiber of $\Xscr$\footnote{More generally, everything that is stated here still holds if $K$ is an arbitrary complete discretely valued field extension of $\QQ_p$ with perfect residue field and if $\Xscr$ is a smooth proper formal scheme over $\Spf(O_K)$ with $X$ its rigid analytic generic fiber (same references). Moreover, some of the results stated below also hold if $\Xscr$ is a proper, flat, and semistable scheme over $O_K$ by work of \v{C}esnavi\v{c}ius and Koshikawa (\cite{CesnaviciusKoshikawa_AInfSemistable}). Conjecturally, for every smooth proper variety $X$ over $K$ there exists such a model $\Xscr$. In this semistable case one has to replace de Rham cohomology by logarithmic de Rham cohomology.}. In this situation there are several $p$-adic cohomology theories attached to it, each of them endowed with some specific interesting additional structure:
\begin{assertionlist}
\item
The integral $p$-adic \'etale cohomology $H^i((X_C)_{\et},\ZZ_p)$ of the generic fiber $X$. Here $C$ is the completion of an algebraic closure $\Kbar$ of $K$. It is a $\ZZ_p$-module with a continuous $\Gamma_K$ action, where $\Gamma_K = \Gal(\Kbar/K)$\footnote{This is a very rich structure compared to the one on the integral $\ell$-adic cohomology $H^i((X_C)_{\et},\ZZ_{\ell})$ ($\ell \ne p$ a prime) that can be identified with $H^i((\Xscr_{\kbar})_{\et},\ZZ_{\ell})$ and whose $\Gamma_K$-action factors through the Galois group of $k$ as explained in Section~\ref{WEIGHTMONODROMY}. This is not true for $H^i((X_C)_{\et},\ZZ_p)$ on which usually the most interesting part of $\Gamma_K$, the pro-$p$-part of the inertia group, acts non-trivially.}.
\item
The de Rham cohomology $H^i_{\dR}(\Xscr/O_K)$ of $\Xscr$. It is an $O_K$-module that is endowed with a filtration, the Hodge filtration (see Section~\ref{RATIONALHODGE}).
\item
The crystalline cohomology $H^i_{\rm crys}(\Xscr_k/W(k))$ defined by Grothendieck that can be thought of as a replacement for the \'etale cohomology $H^i(\Xscr_k, \ZZ_p)$ which formally exists but is usually too small to be of much interest. It is a module over the ring of Witt vectors $W(k)$ equipped with a Frobenius operator which is invertible after inverting $p$.
\end{assertionlist}

Bhatt, Morrow, and Scholze realized in \cite{BMS_IntegralHodge} (see also \cite{BMS_THH} and \cite{Bhatt_Specializing}) that all these cohomology theories are in fact specializations of a single cohomology theory, that was called $A_{\inf}$-cohomology because this cohomology had values in the ring $A_{\inf}$. Recently Bhatt and Scholze generalized and simplified this construction and introduced prismatic cohomology in \cite{BS_Prismatic} (see also Bhatt's lecture \cite{Bhatt_LecturePrismatic}). The next subsections will attempt to give some idea of this construction and its relation to the more classical theory. 


\subsection{Prisms}

Let $R$ be a ring. By considering $\Spec R$ we view $R$ as a geometric object of, say, dimension $d$\footnote{Here the Krull dimension of $\Spec R$ is meant and we do not worry about its various pathologies, such as that $\Spec R$ might be infinite dimensional.}. An already very interesting example is $R = O_K$, where $K$ is some non-archimedean field containing $\QQ_p$, for instance $R = \ZZ_p$ for $K = \QQ_p$. In this case $R$ should be viewed as a $1$-dimensional object.

A prism for $R$ can be thought of as a $(d+1)$-dimensional thickening $\Spec A$ ($A$ a ring) containing $\Spec R$ as a divisor which is also endowed with a ``Frobenius'', i.e., a map $\phi$ that lifts the Frobenius on $A/pA$. In fact, because of technical reasons, it is more convenient not to work with Frobenius lifts $\phi$ but with maps $\delta\colon A \to A$ that measure the difference between the ring endomorphism $\phi$ and the map $x \sends x^p$, i.e. $\phi(a) = a^p + p\delta(a)$. Clearly $\delta$ determines $\phi$. Conversely, $\phi$ determines $\delta$ if $A$ is $p$-torsion free, i.e., the multiplication with $p$ on $A$ is injective.

\begin{definition}\label{DefPrism}
A \emph{prism} is a triple $(A,\delta\colon A \to A ,A \to \Abar)$ satisfying the following conditions.
\begin{definitionlist}
\item
The map $\delta$ satisfies $\delta(0) = \delta(1) = 0$, $\delta(ab) = a ^p\delta(b) + b^p\delta(a) + p\delta(a)\delta(b)$, and $p(\delta(a+b) - \delta(a) - \delta(b)) = a^p + b^p - (a+b)^b$, hence $\delta$ has those properties that ensure that $\phi$, defined by $\phi(a) = a^p + p\delta(a)$, is a ring endomorphism.
\item
The ring map $A \to \Abar$ is surjective, locally on $\Spec A$ its kernel $I$ is generated by a nonzerodivisor $d \in A$, and one has $p \in I + \phi(I)A$.
\item
The ring $A$ is $(p,I)$-adically complete.

Strictly speaking, it has to be derived $(p,I)$-adically complete. By \cite[3.7]{BS_Prismatic}, one can forget about this technical point in case there exists a positive integer $c$ such that $\set{x \in \Abar}{p^nx=0} = \set{x \in \Abar}{p^cx=0}$ for all $n \geq c$. This will be always the case in all examples here and we will assume it from now on\footnote{Such prisms are called \emph{bounded} in \cite{BS_Prismatic}. In other words we assume from now on that all prisms are bounded.}.
\end{definitionlist}
\end{definition}

These conditions on a prism are quite natural. The condition that $A$ is $(p,I)$-adically complete means that one can consider $\Spec A$ simultaneously as a thickening of $A/pA$ and of $R = \Abar = A/I$. Let us consider the condition that $p \in I + \phi(I)A$. By definition, the vanishing scheme $D := \Spec \Abar$ is a divisor on $\Spec A$. The ideal $\phi(I)A$ corresponds to the closed subscheme $V(\phi(I)A) = \phi^{-1}(D)$, where $\phi\colon \Spec A \to \Spec A$ also denotes the map of schemes corresponding to $\phi$. As modulo $p$ the map $\phi$ is the Frobenius which induces the identity on topological spaces, it is clear that $D \cap \phi^{-1}(D)$ contains the characteristic $p$ locus $V(pA) \cap D$ within $D$ (set theoretically). The condition $p \in I + \phi(I)A$ implies that conversely one has $D \cap \phi^{-1}(D) \subseteq V(pA)$.

There is the obvious notion of a map of prisms $(A,\delta,A \to \Abar) \to (B,\delta,B \to \Bbar)$ and one automatically has $\Bbar = B \otimes_A \Abar$ (\cite[3.5]{BS_Prismatic}).

\begin{example}\label{ExamplePrism}
Let $(A,\delta,A \to \Abar)$ be a prism with attached Frobenius lift $\phi$ and denote by $I$ the kernel of $A \to \Abar$.
\begin{assertionlist}
\item\label{ExamplePrism1}
To every perfectoid ring $R$ one can attach the prism $(A_{\inf}(R),\phi,\theta)$ with the notation explained in Subsection~\ref{PERFDRING}. Then $R = \Abar$. Here $\phi$ is the Frobenius on $W(R^{\flat})$ which is bijective. In fact this construction yields an equivalence between the category of perfectoid rings and the category of prisms with bijective Frobenius lift (\cite[3.10]{BS_Prismatic}). In this case the kernel of $\theta$ is a principal ideal.
\item\label{ExamplePrism2}
Let $K$ be a discretely valued complete extension of $\QQ_p$ (e.g., if $K$ is a finite extension of $\QQ_p$), let $\Abar := O_K$ be its ring of integers, choose a prime element $\pi$ of $\Abar$, and let $k = \Abar/\pi \Abar$ be its residue field that we assume to be perfect. The universal property of the ring of Witt vectors $W := W(k)$ (see Subsection~\ref{WITT}) yields an embedding $W \subseteq O_K$. Set $A := W\dlbrack u \drbrack$, the ring of formal power series of $W$ in the variable $u$. Its $\delta$-structure is given by the Frobenius lift which is the Frobenius on the Witt vectors on $W$ and which sends $u$ to $u^p$. Let $A \to \Abar$ be the unique map of $W$-algebras that sends $u$ to $\pi$. Then $(A,\delta,A \to \Abar)$ is a prism and the kernel of $A \to \Abar$ is generated by $E(u) \in A$, where $E$ is an Eisenstein polynomial for $\pi$.
\end{assertionlist}
\end{example}


\subsection{Prismatic cohomology and prismatic site}

From now on, we fix a prism $(A,\delta,A \to \Abar)$ and a smooth\footnote{Everything below may be extended -- although not verbatim -- to not necessarily smooth $p$-adic formal schemes $X$, see \cite[7.2]{BS_Prismatic}.} $p$-adic formal scheme $X$ over $\Abar$. Bhatt and Scholze define in this situation the \emph{prismatic cohomology} $H^i_{\Prism}(X/A)$. It is an $A$-module endowed with an endomorphism $\Phi = \Phi_{H^i}$ that is $\phi$-linear, i.e. $\Phi(m+n) = \Phi(m) + \Phi(n)$ and $\Phi(am) = \phi(a)\Phi(m)$ for all $a \in A$ and $m,n \in H^i_{\Prism}(X/A)$. Moreover, after inverting the ideal $I$ the linearization of $\Phi$ becomes an isomorphism.

In fact, they first define the prismatic site $(X/A)_{\Prism}$ and the structure sheaf $\Oscr_{\Prism}$, and then the prismatic cohomology is defined as the cohomology of that sheaf on the prismatic site.

The \emph{prismatic site} $(X/A)_{\Prism}$ is by definition the category of maps $(A,\delta,A \to \Abar) \to (B,\delta,B\to \Bbar)$ of prisms together with a map $\Spf(\Bbar) \to X$ over $\Abar$ endowed with the flat topology (see~\cite[4.1]{BS_Prismatic} for details). Its structure sheaf $\Oscr_{\Prism}$ sends $(B,\delta,B\to \Bbar, \Spf(\Bbar) \to X)$ to the $A$-algebra $B$. The Frobenius lifts on $B$ define a Frobenius endomorphism $\phi$ on $\Oscr_{\Prism}$. One also defines the sheaf $\overline{\Oscr}_{\Prism}$ of $\Abar$-algebras, that sends $(B,\dots)$ to $\Bbar$. Both are indeed sheaves by \cite[3.12]{BS_Prismatic}.



\subsection{Hodge Tate comparison}

The Hodge-Tate comparison relates prismatic cohomology and sheaves of differentials. It is the main tool for controlling prismatic cohomology.

\begin{theorem}(\cite[4.10]{BS_Prismatic})\label{HodgeTateComp}
Let $X = \Spf(S)$ be an affine\footnote{This is only a technical assumption to avoid introducing more notation.} formal scheme. Then there are functorial isomorphisms of $S$-modules
\[
(\Omega^i_{S/\Abar})^{\wedge}\{i\} \cong H^i((X/A)_{\Prism},\overline{\Oscr}_{\Prism}).
\]
\end{theorem}

Here $(\cdot)^{\wedge}$ denotes $p$-adic completion and $(\cdot)\{i\}$ denotes a twist that can be ignored if $I$ is a principal ideal (which is usually the case) and one has chosen a generator of $I$. This result implies that the formation of prismatic cohomology is compatible with base change (\cite[4.11]{BS_Prismatic}) and commutes with \'etale localization (\cite[4.19]{BS_Prismatic})\footnote{But note that both results are shown in their affine version in \cite{BS_Prismatic} without using the Hodge-Tate comparison and are used as ingredients to prove the Theorem~\ref{HodgeTateComp}.}. It also gives a description of the cotangent complex of $X$ over $A$ in terms of prismatic cohomology (\cite[4.14]{BS_Prismatic}, see also \cite[3.2.1]{ALB_Dieudonne}).


\subsection{Integral $p$-adic cohomology theories as specializations of prismatic cohomology}

The main result of \cite{BS_Prismatic}, extending results of \cite{BMS_IntegralHodge} to prismatic cohomology, is the fact that crystalline cohomology, de Rham cohomology and \'etale cohomology can all be obtained by specializing the prismatic cohomology.

\begin{theorem}(\cite[1.8]{BS_Prismatic})\label{PrismaticSpecialize}
Let $(A,\delta,A \to \Abar)$ be a prism, $\phi$ the corresponding Frobenius lift on $A$, $I$ the kernel of $A \to \Abar$, and let $X$ be a smooth $p$-adic formal scheme over $\Abar$.
\begin{assertionlist}
\item\footnote{This crystalline comparison theorem is proved in \cite{BS_Prismatic} before the Hodge-Tate comparison theorem and it is used as an ingredient of its proof.}
If $I = (p)$, then the Frobenius pullback of prismatic cohomology yields crystalline cohomology.
\item
``Specializing''\footnote{To make ``specializing'' precise one should work with the complex $R\Gamma_{\Prism}(X/A)$ in the derived category of $A$-modules whose cohomology prismatic cohomology and consider derived base change.} prismatic cohomology via the composition $A \ltoover{\phi} A \to \Abar$ one obtains the de Rham cohomology of $X$ over $\Abar$.
\item
Assume that $\phi$ is an isomorphism. Then the prism is automatically as in the situation of Example~\ref{ExamplePrism}~\ref{ExamplePrism1}. In particular, $\Abar$ is perfectoid. Let $X_{\eta}$ be the adic generic fiber of $X$ over $\Abar[1/p]$. Then for any $n \geq 1$ the \'etale cohomology of $X_{\eta}$ with values in $\ZZ/p^n\ZZ$ is obtained from ``specializing'' the prismatic cohomology modulo $p^n$, inverting $I$, and taking $\Phi$-invariants.
\end{assertionlist}
\end{theorem}

Suppose that one is in the situation of Example~\ref{ExamplePrism}~\ref{ExamplePrism2} and that $X$ is a smooth proper formal scheme over $O_K$, then the prismatic cohomology of $X/A$ are finitely generated $A$-modules endowed with a Frobenius map that becomes an isomorphism after inverting the Eisenstein polynomial $E$. Such modules have also been defined by Breuil and studied further by Kisin. They are accordingly called \emph{Breuil-Kisin modules} and play an important role in the study of $p$-adic Galois representations. In special cases (e.g., if $X$ is an abelian scheme over $O_K$, or if the cohomology degree is smaller than $p-1$), there were previously other constructions of Breuil-Kisin modules attached to smooth proper formal schemes over $O_K$ often using the \'etale $p$-adic cohomology of the generic fiber as an ingredient, for instance in \cite{CaisLiu_BKModules}. The results of Bhatt and Scholze show that in complete generality these constructions are ``shadows'' of a deeper structure, the prismatic cohomology.


\subsection{Further applications of prismatic cohomology}

The formalism of prismatic cohomology has several further applications. We mention only three of them briefly.

\subsubsection{Perfectoidization of integral algebras}\label{PERFDIZATION}

One application is the existence of a canonical perfectoidization of rings that are integral over a perfectoid ring.

\begin{theorem}(\cite[1.16]{BS_Prismatic}\label{Perfectoidization}
Let $R$ be a perfectoid ring and let $S$ be an integral $R$-algebra. There is a map of rings $S \to S_{\rm perfd}$ with $S_{\rm perfd}$ perfectoid such that any map of rings from $S$ to a perfectoid ring factors uniquely through $S_{\rm perfd}$.
\end{theorem}

In fact, Bhatt and Scholze define in \cite[8.2]{BS_Prismatic} a perfectoidization for every algebra\footnote{even simplicial algebra} over a perfectoid ring using derived prismatic cohomology. They obtain a (derived) $p$-complete $\EE_{\infty}$-ring over $S$ and show that if this is a usual ring (i.e., concentrated in degree $0$), then $S_{\rm perfd}$ is perfectoid (\cite[8.13]{BS_Prismatic}). This is the case if $R \to S$ is the $p$-completion of an integral map (\cite[10.11]{BS_Prismatic}).

\subsubsection{Nygaard filtration}

Bhatt and Scholze define on the Frobenius twist of prismatic cohomology a filtration, the Nygaard filtration, and describe its graded pieces (see \cite[15.3]{BS_Prismatic} for details). It lifts the conjugate filtration defined by the Hodge-Tate comparison.

\subsubsection{Prismatic Dieudonn\'e theory}

In \cite{ALB_Dieudonne} Ansch\"utz and Le Bras show that over a very general class\footnote{the class rings that are are flat over $\ZZ_p$ or over some $\ZZ/p^n\ZZ$ and that are quasi-syntomic (e.g., any noetherian $p$-complete locally complete intersection ring is quasi-syntomic)} of rings $R$ one can classify $p$-divisible groups by their filtered Dieudonn\'e modules. This generalizes all previously known statements of Dieudonn\'e theory. The precise formulation and the proof uses heavily the machinery of prisms and prismatic cohomology.


\part{Langlands program}

Let $F$ be a finite extension of $\QQ$ and let $\AA_F$ be its ring of adeles. As before, fix a prime $p$. Choose an isomorphism of fields $\QQbar_p \cong \CC$. As explained above in Subsection~\ref{INTROLANGLANDS}, the global Langlands conjecture makes the following prediction.

\begin{conjecture}\label{GlobalLanglands}
There is a bijection\footnote{This bijection has to preserve certain invariants. More precisely, it should match eigenvalues of Frobenius elements with Satake parameters.} of isomorphism classes of certain\footnote{irreducible continuous, almost everywhere unramified and de Rham at places dividing $p$} $n$-dimensional representations of $\Gal(\Fbar/F)$ over $\QQbar_p$ and of certain\footnote{$L$-algebraic, cuspidal, automorphic} representations of $\GL_n(\AA_F)$ over $\CC$.
\end{conjecture}

There is also a local version if $F$ is a finite extension of $\QQ_p$ which is in fact a theorem by work of Henniart and Harris-Taylor.

\begin{theorem}\label{LocalLanglands}
There is a bijection\footnote{compatible with twists, central characters, duals, $L$- and $\epsilon$-factors in pairs} of isomorphism classes of certain\footnote{supercuspidal, irreducible, smooth} representations of $\GL_n(F)$ over $\CC$ and of irreducible $n$-dimensional continuous representations of the Weil group\footnote{The Weil group is a variant of the Galois group that allows Frobenius elements in the Galois group of $F$ (as explained in Subsection~\ref{WEIGHTMONODROMY}) to act arbitrarily on continuous representations over the complex numbers.} $W_F$ of $F$ over $\CC$.
\end{theorem}

The global and the local Langlands correspondence have generalizations by replacing $\GL_n$ by an arbitrary reductive group $G$ and reformulating the notion of an $n$-dimensional Galois/Weil representation accordingly. In full generality this is still conjectural even in the local case. Moreover, the principle of Langlands functoriality predicts roughly that these bijections should be ``functorial in $G$''.

As mentioned above, the analogues of these conjectures over global and local fields of characteristic $p$ are known now by the work of many people (Drinfeld, Laumon, Rapoport, Stuhler, L.~Lafforgue, V.~Lafforgue, and others). Moreover, the local Langlands conjecture is also known for some classical groups, among them $\GL_n$ (by work of Henniart, Harris, and Taylor). The global Langlands conjecture is still wide open in general except for $G = \GL_1$ where it is essentially equivalent to global class field theory.

\subsubsection*{Moduli spaces of Shtukas and Shimura varieties}

The main tool to establish a Langlands correspondence is the construction of geometric objects that have so many symmetries such that their cohomology carries a Galois action and an automorphic action. And then one hopes to prove that this cohomology realizes the Langlands correspondence.

To prove the global Langlands conjecture for global fields $F$ of characteristic $p$, one first observes that $F$ corresponds by a standard result in Algebraic Geometry to a curve $X$ defined over a finite field (e.g., \cite[15.22]{GW}). Now one uses the moduli spaces $G$-${\rm Sht}_{X,I}$ of so-called ``$G$-shtukas with legs over $X$'', where the legs are parametrized by a finite set $I$. These are moduli spaces of $G$-bundles $\Escr$ on the curve $X$ together with a Frobenius linear automorphism of $\Escr$ defined outside a family $(x_i)_{i\in I}$ of points on $X$. These moduli spaces exist as inductive limits of algebraic stacks for an arbitrary reductive group $G$. Attaching to a shtuka its legs defines a morphism $G$-${\rm Sht}_{X,I}$ to $I$-fold self product $X^I$ of $X$.  

For global fields $F$ of characteristic $0$ there is unfortunately no analogue for $G$-${\rm Sht}_{X,I}$ in general. Here one of the best available geometric objects are Shimura varieties. But they have several drawbacks. They can be attached only to a particular class of reductive groups. For instance, there is no Shimura variety for $\GL_n$ over $\QQ$ if $n \geq 3$. Moreover, Shimura varieties miss the important structure of legs as defined for shtukas. Here $X$ would be $\Spec O_F$ (or some ``Arakelov completion'' of it taking into account the archimedean places) but it is totally unclear how to define $(\Spec O_F)^I$, which ideally should be the $I$-fold fiber product of $\Spec O_F$ over the non-existing ``field with one element''.

Nevertheless, Shimura varieties have turned out to be a central tool in Arithmetic Geometry. For instance the theory of modular curves, which are the special case of Shimura varieties attached to $\GL_2$, has been fundamental in the Taylor-Wiles proof of Fermat's Last Theorem. Shimura varieties, attached to certain unitary groups, have also been the central tool in the proof of Henniart and Harris-Taylor of the local Langlands conjecture for $\GL_n$ over $p$-adic fields.

Shimura varieties are quite complicated objects that combine algebraic geometric, arithmetic, analytic, and representation theoretic aspects. Here the following properties will be important.
\begin{assertionlist}
\item
A Shimura variety is attached to a datum given by a triple $(G,X,K)$, where $G$ is a reductive group over $\QQ$, $X$ is a complex manifold\footnote{More precisely, a finite disjoint union of hermitian symmetric domains} of the form $X = G(\RR)/K_{\infty}$ for a maximal compact subgroup $K_{\infty}$ of the real Lie group $G(\RR)$, and $K$ is an open compact subgroup of $G(\AA_f)$ where $\AA_f$ is the topological ring of finite adeles of $\QQ$, called the \emph{level structure}.

The complex manifold $X$ admits a holomorphic embedding, the Borel embedding, into a (smooth projective) flag variety $\Fcal\ell_{G,X}$ which is of the form $G_{\CC}/P_X$ for some parabolic subgroup $P_X$ of $G_{\CC}$. 
\item
The Shimura variety $S_K := S_K(G,X)$ is a smooth quasi-projective variety\footnote{At least if $K$ is sufficiently small. Otherwise one should work with Deligne-Mumford stacks if one does not want limit oneself to coarse moduli spaces.} defined over a subfield $E$ of $\CC$, the \emph{reflex field} of $(G,X)$. The reflex field is a finite extension of $\QQ$ which depends only on $(G,X)$. The complex manifold $S_K(\CC)$ can be described as
\[
S_K(\CC) = G(\QQ)\backslash (X \times G(\AA_f)/K).
\]
\item
The Shimura variety has a natural minimal compactification $S^*_K$, the Baily-Borel-Satake compactification, which is a normal projective variety.
\end{assertionlist}

Examples of Shimura varieties are modular curves, say with full level-$N$ for some $N \geq 3$. Here one chooses $G = \GL_2$,
\[
X = \set{\tau \in \CC}{{\rm Im}(\tau) \ne 0} = \PP^1(\CC) \setminus \PP^1(\RR),
\]
and $K = \set{g \in \GL_2(\widehat{\ZZ})}{g \equiv 1\pmod{N}}$, where $\widehat{\ZZ} = \lim_n \ZZ/n\ZZ = \prod_p\ZZ_p$ is the pro-finite completion of $\ZZ$. These modular curves are moduli spaces of elliptic curves $E$ with a rigidification of their $N$-torsion. In this case the flag variety $\Fcal\ell_{G,X}$ is the projective line $\PP^1_{\CC}$ and the Borel embedding is the inclusion $X \mono \PP^1(\CC)$.

More generally, Shimura varieties attached to $G = \GSp_{2g}$ and
\[
X = \set{\tau \in M_g^{\rm sym}(\CC)}{\text{${\rm Im}(\tau)$ positive or negative definite}}
\]
are moduli spaces of polarized abelian varieties of dimension $g$. These are sometimes called Siegel Shimura varieties. In this case the flag variety $\Fcal\ell_{G,X}$ parametrizes maximal totally isotropic subspaces in a symplectic space of dimension $2g$.


\section{Perfectoidization of Shimura varieties and the Hodge-Tate period map}\label{PERFDSHIMURA}

Scholze showed that there is a natural way to make Shimura varieties of Hodge type\footnote{Scholze conjectures in \cite{Scholze_PerfdShimura} that similar results should hold for arbitrary Shimura varieties. Shen generalized Scholze's results to the large class of Shimura varieties of abelian type in \cite{Shen_PerfectoidAbelianType}.} (which means roughly that the Shimura variety can be embedded into a Siegel Shimura variety) perfectoid. 

For this fix a prime number $p$, a place $\pfr$ of the reflex field $E$ over $p$, and let $E_{\pfr}$ be the completion with respect to the $\pfr$-adic absolute valuation on $E$. By base change we consider $S_K$ as a smooth adic space over $E_{\pfr}$. Fix a sufficiently small compact open subgroup $K^p$ of $G(\AA^p_f)$, where $\AA^p_f$ is the ring of finite adeles of $\QQ$ with trivial $p$-th component.

\begin{theorem}(\cite[2.1.2, 2.1.3]{CS_GenericI})\label{PerfectoidShimura}
The limit\footnote{in a suitable sense, see \cite[4.1]{SW_ModuliPDiv} for a thorough discussion} of $S_{K^pK_p}$, where $K_p$ runs through the open compact subgroups of $G(\QQ_p)$ is a perfectoid space. Moreover, there is a \emph{Hodge-Tate period map}
\[
\pi_{\rm HT}\colon S_{K^p} \to \Fcal\ell^{\rm opp}_{(G,X)},
\]
where $\Fcal\ell^{\rm opp}_{(G,X)}$ is the flag variety opposite to $\Fcal\ell_{(G,X)}$, viewed as an adic space. The Hodge-Tate period map is equivariant for the $G(\QQ_p)$-action and for the prime-to-$p$ Hecke operators for varying $K^p$.
\end{theorem}

These constructions extend to (a variant of) the Baily-Borel-Satake compactification (\cite[5.1]{Scholze_PerfdShimura}).

One can think of the Hodge-Tate period map as a $p$-adic analogue of the Borel embedding.
In the case of the modular curve, $S_{K^p}$ can be described as a moduli space of elliptic curves $E$ together with a trivialization of its Tate module $T_p(E) = H^1(E,\ZZ_p)\vdual$ as $\ZZ_p$-module. In this case $\Fcal\ell^{\rm opp}_{(G,X)}$ is again the projective line $\PP^1$, considered as adic space. If $C$ is an algebraically closed non-archimedean extension of $E_{\pfr}$, one can describe the Hodge-Tate period map on $C$-valued points as the map that sends an elliptic curve $E$ together with a trivialization $T_p(E) \iso \ZZ_p^2$ to the dual of the Hodge-Tate filtration in
\[
C^2 \cong T_p(E) \otimes_{\ZZ_p} C = (H^1(E,\ZZ_p) \otimes_{\ZZ_p} C)\vdual
\]
given by the degenerating Hodge-Tate spectral sequence (Theorem~\ref{HodgeTateSS}). Moreover, one has a decomposition of the adic projective line
\begin{equation}\label{EqModCurveNewton}
(\PP^1)^{\rm ad} = \PP^1(\QQ_p) \sqcup \Omega,
\end{equation}
where $\Omega = (\PP^1)^{\rm ad} \setminus \PP^1(\QQ_p)$ denotes Drinfeld's half plane. Then the preimage of $\Omega$ under the $\pi_{\rm HT}$ is the open tube of the supersingular locus of the special fiber\footnote{More precisely, it consists of those points whose rank-$1$-generizations specialize to a supersingular point of the special fiber.}, and the preimage of $\PP^1(\QQ_p)$ is (the closure of) the tube of the ordinary locus.

A similar description of the Hodge-Tate period map holds more generally for Shimura varieties that parametrize abelian varieties with additional structure, e.g. in the Siegel case. In fact, there is the following theorem of Scholze and Weinstein that can be viewed as a $p$-adic analogue of Riemann's classification of complex abelian varieties via their integral singular homology and their Hodge filtration.

\begin{theorem}(\cite[B]{SW_ModuliPDiv})\label{PDiv}
Let $C$ be an algebraically closed non-archimedean extension of $\QQ_p$. Attaching to a $p$-divisible group $G$ over $O_C$ its Tate module $T$ and $W = \Lie(G) \otimes_{O_C} C$ yields an equivalence between the category of $p$-divisible groups over $O_C$ and the category of pairs $(T,W)$, where $T$ is a free $\ZZ_p$-module and $W \subseteq T \otimes C$ is a $C$-subvectorspace.
\end{theorem}

This shows that one can view the adic flag variety (or rather its $C$-valued points) as a parameter space for $p$-divisible groups with additional structure over $O_C$. The decomposition \eqref{EqModCurveNewton} can then be reinterpreted as the stratification according to the Newton polygon of the special fiber of the $p$-divisible group. And the Hodge-Tate period map, restricted to the locus of elliptic curves of good reduction, can be viewed on $C$-valued points as the map that sends an elliptic curve over $C$ to the $p$-divisible group of its model over $O_C$. 

All of this can be extended to Shimura varieties of Hodge type, see \cite{CS_GenericI}.


\section{Construction of Galois representations}

\subsection{Construction of torsion Galois representations}

Scholze proves in \cite{Scholze_Torsion} a mod-$p$-version of the surjectivity of the bijection in the global Langlands conjecture~\ref{GlobalLanglands} for special number fields. Let $F$ be a totally real field or a field with complex multiplication, let $n \geq 1$ be an integer. A general strategy to relate automorphic forms and Galois representations also $p$-adically is to study the singular cohomology of the locally symmetric spaces for $\GL_n$ over $F$. More precisely, let $K$ be a sufficiently small open compact subgroup $K \subset \GL_n(\AA_{F,f})$, where $\AA_{F,f}$ denotes the ring of finite adeles of $F$. Consider the locally symmetric space
\[
S_K = \GL_n(F)\backslash(X \times \GL_n(\AA_{F,f})/K),
\]
where $X = \GL_n(F \otimes_{\QQ} \RR)/\RR_{>0}K_{\infty}$ is the symmetric space for $\GL_n(F \otimes_{\QQ} \RR)$ and where $K_{\infty} \subset \GL_n(F \otimes_{\QQ} \RR)$ is a maximal compact subgroup. Note that $X$ in general does not carry a complex structure and $S_K$ is not a Shimura variety. The singular cohomology spaces $H^i(S_K,\CC)$ carry an action by an algebra of Hecke operators and all Hecke eigenvalues in $H^i(S_K,\CC)$ come up to a twist from certain automorphic representations. Hence the following theorem of Scholze can be viewed as a modulo $p$ variant of the global Langlands conjecture for totally real or CM fields $F$.

\begin{theorem}(\cite[1.0.3]{Scholze_Torsion})\label{TorsionLanglands}
For any system of Hecke eigenvalues appearing in $H^i(S_K,\FF_p)$ there is a continuous semi-simple representation $\Gal(\Fbar/F) \to \GL_n(\FFbar_p)$ characterized by the property that for all but finitely many ``ramified'' places $v$ of $F$, the characteristic polynomial of the Frobenius at $v$ is described in terms of the Hecke eigenvalues at $v$.
\end{theorem}

Scholze proves also a version with coefficients in $\ZZ/p^m\ZZ$. Passing to the limit over $m$ he obtains as a consequence the existence of $n$-dimensional Galois representations of $\QQbar_p \cong \CC$ attached to certain representations\footnote{cuspidal automorphic representations $\pi$ such that $\pi_{\infty}$ is regular $L$-algebraic} of $\GL_n(\AA_F)$, which had been previously shown by Harris-Lan-Taylor-Thorne. But it is important to note that the dimension of $H^i(S_K,\FF_p)$ is in general much larger than that of $H^i(S_K,\CC)$ and most systems of Hecke eigenvalues cannot be lifted to characteristic $0$.

The key automorphic result in Scholze's proof is that all Hecke eigenvalues in the compactly supported completed cohomology of a Shimura variety of Hodge type can be $p$-adically interpolated by Hecke eigenvalues coming from classical cusp forms.
Its proof proceeds by realizing the cohomology of $S_K$ as the boundary contribution of a Shimura variety of Hodge type and then using its perfectoidization and the Hodge-Tate period map (Section~\ref{PERFDSHIMURA}).


\subsection{Potential automorphy over CM fields}

Theorem~\ref{TorsionLanglands} is an essential tool used by Allen, Calegari, Caraiani, Gee, Helm, Le Hung, Newton, Scholze, Taylor, and Thorne in \cite{ACCGHLNSTT}. Here the first unconditional modularity lifting theorems for $n$-dimensional ($n$ arbitrary!) regular Galois representations without any self-duality conditions is shown.

As an application they obtain the potential modularity of arbitrary elliptic curves (and their symmetric powers) over CM fields. In particular the Sato-Tate conjecture holds for such elliptic curves.

As another application they obtain the Ramanujan conjecture\footnote{The (generalized) Ramanujan conjecture predicts that all local components of cuspidal automorphic representations are tempered.} for certain cuspidal automorphic representations of $\GL_2(\AA_F)$ for all CM fields $F$.

The arguments in \cite{ACCGHLNSTT} crucially exploit the work of Caraiani and Scholze in \cite{CS_GenericII} on the cohomology of non-compact Shimura varieties associated to unitary groups.


\section{Local Langlands in mixed characteristic}\label{LOCALLANGLANDS}

\subsection{Classical local Langlands correspondence for $\GL_n$ over local fields}

The local Langlands correspondence~\ref{LocalLanglands} was proved by Harris-Taylor and by Henniart at the end of the 90s. Both proofs used in an essential way Henniart's numerical local Langlands conjecture which is based on a complicated reduction modulo $p$. It implies that it suffices to construct any (not necessarily bijective) map from the set of (supercuspidal irreducible smooth) representations of $\GL_n(F)$ to the set of $n$-dimensional Weil group representations satisfying all required compatibilities.

Scholze gave a new and much simpler proof that completely bypasses the numerical local Langlands conjecture in his master thesis (\cite{Scholze_LLC}). It is based on a geometric argument via the nearby cycles sheaves of certain moduli spaces of $p$-divisible groups.


\subsection{The Fargues-Fontaine curve}\label{FFCURVE}

The classical local Langlands correspondence is just a highly non-trivial bijection of sets. This is much less satisfactory than the versions of local and global Langlands correspondence in characteristic $p$ that yield a geometric interpretation and a categorification of the Langlands correspondence. 

In recent years much progress was made to obtain a similar (still conjectural) geometrization also over $p$-adic fields. Recall that in characteristic $p$ the Langlands correspondence was obtained via the study of a moduli space of shtukas which are given by $G$-bundles (and additional data) on a curve. We refer to \cite{Fargues_GeometrizationLLC} for details of the fascinating development. Here we will explain only the definition of the $p$-adic curve involved in this as it will be also essential in the theory of $p$-adic shtukas and local Shimura varieties explained below. It is the Fargues-Fontaine curve.

\subsubsection*{Definition of the Fargues-Fontaine curve}
In fact, we recall from \cite[3.3]{CS_GenericI} how to attach to every perfectoid space $S$ in characteristic $p$ an adic space $X_S$, called the \emph{Fargues-Fontaine curve attached to $S$}. First assume that $S = \Spa(R,R^+)$ is affinoid. Then $R^+$ is a perfect ring. The adic spectrum $\Spa (W(R^+),W(R^+))$ should be thought of having two coordinates\footnote{But there is no good rigorous notion of dimension in this setting!} given by the functions $p$ and $[\varpi]$ in $W(R^+)$, where $\varpi \in R^+$ is a pseudo-uniformizer. Let $Y_S$ be the locus in $\Spa (W(R^+),W(R^+))$ where the product of coordinate functions $p[\varpi]$ is non-zero. The Frobenius on $W(R^+)$ induces an automorphism $\varphi$ on $Y_S$, and the Fargues-Fontaine curve attached to $S$ is defined as the quotient of $Y_S$ be this action
\[
X_S := Y_S/\varphi^{\ZZ}.
\]
This can be shown to be an adic space again. One then checks that the functor $S = \Spa(R,R^+) \sends Y_S$ glues to a functor $S \sends Y_S$ for $S$ a perfectoid spaces in characteristic $p$ and one obtains in general $X_S := Y_S/\varphi^{\ZZ}$ as an adic space. It is an analytic adic space over $\QQ_p$.

If $S = \Spa (C,C^{\circ})$, where $C$ is an algebraically closed non-archimedean field of characteristic $p$, then $X_S = X_C$ should be thought of as a curve. In fact, originally Fargues and Fontaine defined the schematic Fargues-Fontaine curve $X^{\rm sch}_C$ which is a noetherian regular $1$-dimensional scheme over $\QQ_p$ (but not of finite type over any field). The adic Fargues-Fontaine curve $X_C$ can be thought of as an analytification of $X_C^{\rm sch}$.

\subsubsection*{Untilts and distinguished divisors}

If $(S^{\sharp},\iota)$ is an untilt of $S$ over $\Spa \QQ_p$, then there exists a closed immersion $i_{S^{\sharp}}\colon S^{\sharp} \to X_S$ which is obtained by composing the projection $Y_S \to X_S$ with $\theta\colon S^{\sharp} \to Y_S$ which is locally given by $W(R^+) \to R^{\sharp+}$ \eqref{EqDefTheta}. The closed immersion $\theta\colon S^{\sharp} \to Y_S$ identifies $S^{\sharp}$ with a \emph{distinguished Cartier divisor}, i.e., $S^{\sharp}$ is locally on $Y_S$ the vanishing locus of a nonzerodivisor $\xi$ in $W(R^+)$ such that $\varphi(\xi) - \xi^p = pu$ for a unit $u\in W(R^+)^{\times}$ (\cite[3.11]{BMS_IntegralHodge}). In fact, this construction induces a bijection between the set of isomorphism classes of untilts of $S$ and the set of distinguished Cartier divisor in $Y_S$ (\cite[1.22]{Fargues_GeometrizationLLC}). 

\subsubsection*{The diamond formula for the Fargues-Fontaine curve}

One should view $Y_S$ as a product ``$S \times \Spa\QQ_p$''. Literally this does not make sense because $S$ lives in characteristic $p$ and $\QQ_p$ is a field of characteristic $0$. But by Scholze's theory of diamonds (Section~\ref{DIAMOND}) we can view objects over $\QQ_p$ as objects in characteristic $p$ and then there are the following isomorphisms of diamonds
\[
Y_S^{\diamond} = S \times \Spd(\QQ_p), \qquad X_S^{\diamond} = (S \times \Spd(\QQ_p))/\varphi^{\ZZ}.
\]
where $\varphi$ on the right hand side is given by the product of the Frobenius on $S$ and the identity on $\Spd(\QQ_p)$.


\subsection{The geometrization of the local Langlands correspondence}\label{GEOMLLC}

Let $G$ be a reductive group over a $p$-adic field $F$. We now state the geometrization of the local Langlands conjecture due to Fargues. The reference is \cite{Fargues_GeometrizationLLC}.

To simplify the notation we assume that $F = \QQ_p$ and that $G$ is quasi-split (i.e., contains a Borel subgroup defined over $\QQ_p$). One can now define the notion of a $G$-bundle on the Fargues-Fontaine curve. For instance, if $G = \GL_n$, then a $\GL_n$-bundle is simply a rank $n$ vector bundle. Similarly, as in the characteristic $p$ variant of the global Langlands correspondence one considers the stack of $G$-bundles on the Fargues-Fontaine curve. More precisely, one attaches to each perfectoid space $S$ the groupoid $\Bun_G(S)$ of $G$-bundles on the Fargues-Fontaine curve $X_S$. This is a stack for the $v$-topology on $\Perf$.

One can construct $G$-bundles on $X_S$ as follows. Let $\QQbreve_p$ be the completion of the maximal unramified extension of $\QQ_p$. It is the field of fractions of $W(\FFbar_p)$. Denote by $\sigma$ its Frobenius and by $B(G)$ the quotient of $G(\QQbreve_p)$ by $\sigma$-conjugacy, i.e., by the equivalence relation $b \sim gb\sigma(b)^{-1}$ for $g \in G(\QQbreve_p)$. Given $b \in G(\QQbreve_p)$ one attaches a $G$-bundle $\Escr_b := (G_{\QQbreve_p} \times_{\QQbreve_p} Y_S)/\varphi^{\ZZ}$, where $\varphi$ acts diagonally by the Witt vector Frobenius on $Y_S$ and by $g \sends b\sigma(g)$ on $G_{\QQbreve_p}$. For instance $\Escr_1$ is the trivial bundle. This construction induces a map
\begin{equation}\label{EqGBundle}
B(G) \lto \{\text{$G$-bundles on $X_S$}\}/\cong.
\end{equation}
If $S = \Spa(C,C^+)$ for an algebraically closed field, then this map is bijective by \cite{Fargues_GBundle}. This implies that the underlying set of $\Bun_G$ can be identified with $B(G)$.

To every $b \in G_{\QQbreve_p}$ one attaches a reductive group $J_b$ over $\QQ_p$ with
\[
J_b(\QQ_p) = \set{g \in G_{\QQbreve_p}}{gb\sigma{g}^{-1} = b}.
\]
The isomorphism class of $J_b$ depends only on the $\sigma$-conjugacy class of $b$ in $B(G)$. An element in $B(G)$ is called \emph{basic} if $J_b$ is an inner form of $G$. For instance $b = 1$ is basic with $J_1 = G$. For basic $b$, the automorphism group of the $G$-bundle $\Escr_b$ is the locally profinite group $J_b(\QQ_p)$ (for general $b$, $J_b(\QQ_p)$ is only a subgroup of the automorphism group). By the general formalism of residual gerbes, this implies that if $\Fscr$ is a $\QQbar_{\ell}$-sheaf on $\Bun_G$, then its stalk at a basic point $b$ is a $\QQbar_{\ell}$-vector space endowed with an action of $J_b(\QQ_p)$.

Now Fargues conjectures in \cite{Fargues_GeometrizationLLC} that the local Langlands correspondence for arbitrary reductive groups can be realized geometrically as follows. There is a functor of the groupoid of representations of the Weil group with values in the $L$-group of $G$\footnote{more precisely, a discrete Langlands parameter of $G$} to the groupoid of certain (perverse) $\QQbar_{\ell}$-sheaves on $\Bun_G$ endowed with an action of the Weil group $W_F$. The stalks of this sheaves at basic points $b$ of $\Bun_G$ are then endowed with actions of $W_F$ and of $J_b(\QQ_p)$. Then the conjecture says that they define a local Langlands correspondence for the inner form $J_b$.


\section{Local Shimura varieties}

The main reference for this section is \cite{SW_Berkeley}.

As Shimura varieties can only be attached to special reductive groups, one might hope to find a replacement that one can define for arbitrary reductive groups over local or even global fields of characteristic $0$. Over local (non-archimedean) fields, Rapoport and Viehmann conjectured in \cite{RV_TowardsLocal} the existence of \emph{local Shimura varieties}. These are attached to arbitrary reductive groups $G$ over $\QQ_p$ and to a cocharacter\footnote{The cocharacter is assumed to be minuscule in \cite{RV_TowardsLocal} but for some of the constructions explained below this hypothesis is not necessary.} $\mu$ of $G$, which corresponds roughly to the datum of the hermitian symmetric domain $X$ for a Shimura variety. Moreover, one fixes an element $b \in G(\QQbreve_p)$ satisfying some group theoretic compatibility with $\mu$. By a theorem of Rapoport this compatibility is equivalent to the non-emptiness of the local Shimura varieties defined below if $\mu$ is minuscule. The field of definition of the conjugacy class of $\mu$ is denoted by $E$. It is a finite extension of $\QQ_p$. We denote by $\Ebreve$ the completion of the maximal unramified extension of $E$.

\subsection{Moduli spaces $p$-adic shtukas}

Scholze and Weinstein attach in \cite[XXIII]{SW_Berkeley} to these data the \emph{moduli space of Shtukas (with one leg)} $\Sht_{G,b,\mu}$. It is the sheaf that attaches to every perfectoid space $S$ in characteristic $p$ the set of pairs $((S^{\sharp},\iota), \alpha)$. Here $(S^{\sharp},\iota)$ is an untilt of $S$. It defines in the Fargues-Fontaine curve $X_S$ a divisor (Subsection~\ref{FFCURVE}) that we simply denote by $S^{\sharp}$. As explained in the previous section, the element $b$ defines a $G$-bundle $\Escr_b$ on the Fargues-Fontaine curve. Then $\alpha$ is an isomorphism of $G$-bundles
\[
\alpha\colon \Escr_1\rstr{X_S \setminus S^{\sharp}} \liso \Escr_b\rstr{X_S \setminus S^{\sharp}}
\]
which is meromorphic along $S^{\sharp}$ and its singularity along $S^{\sharp}$ is bounded by $\mu$.

Composition of $\alpha$ with automorphisms of $\Escr_1$ and of $\Escr_b$ defines a continuous action of the locally profinite group $G(\QQ_p) \times J_b(\QQ_p)$ on $\Sht_{G,\mu,b}$. If $b$ is basic, the situation is symmetric\footnote{Strictly speaking, this is only true if one formulates the entire theory for general not necessarily quasi-split reductive groups $G$. If $G$ is quasi-split, then $J_b$ is only quasi-split if $J_b = G$, which is usually not the case.} in $G$ and $J_b$. One obtains a duality between moduli spaces of shtukas, namely a $G(\QQ_p) \times J_b(\QQ_p)$-equivariant isomorphism
\[
\Sht_{G,b,\mu} \cong \Sht_{\breve{G}, \breve{b},\breve{\mu}}
\]
of diamonds over $\Spd \Ebreve$ (\cite[23.3.2]{SW_Berkeley}). Here $\breve{G} = J_b$, $\breve{b} = b^{-1}$, and $\breve{\mu} = \mu^{-1}$ under the identification $\breve{G}_{\QQbreve_p} = G_{\QQbreve_p}$. This generalizes work of Faltings and Fargues for $G = \GL_n$ and certain $\mu$. Furthermore, this duality isomorphism exchanges the Hodge-Tate period map with the de Rham period map.

If $K \subseteq G(\QQ_p)$ is an open compact subgroup, we set $\Sht_{G,\mu,b,K} := \Sht_{G,\mu,b}/K$. Then $\Sht_{G,\mu,b} = \lim_K \Sht_{G,\mu,b,K}$. Moreover, one has the following theorem.

\begin{theorem}(\cite[23.3.4]{SW_Berkeley})\label{ShtukaDiamond}
The sheaf $\Sht_{G,\mu,b,K}$ on $\Perf$ is a diamond.
\end{theorem}

To show the theorem, one constructs a Grothendieck-Messing period map from $\Sht_{G,\mu,b,K}$ to a mixed characteristic variant of an affine Grassmannian ${\rm Gr}_{G,\leq\mu}$ and shows that it is \'etale. As ${\rm Gr}_{G,\leq\mu}$ is a diamond (\cite[20.2.3]{SW_Berkeley}), one deduces that $\Sht_{G,\mu,b,K}$ is a diamond.


\subsection{Local Shimura varieties}

From now on we assume that $\mu$ is minuscule, i.e. $\langle \mu,\alpha\rangle \in \{-1,0,1\}$ for every root $\alpha$ in the based root system of $G_{\kbar}$. In this case, the mixed characteristic variant of an affine Grassmannian ${\rm Gr}_{G,\leq\mu}$ is the diamond of a flag variety $\Fcal\ell_{G,\mu}$ attached to $G$ and $\mu$. From the existence of the Grothendieck Messing period map one deduces that there exists a unique smooth rigid analytic space $\Mscr_{G,\mu,b,K}$ such that
\[
(\Mscr_{G,\mu,b,K})^{\diamond} = \Sht_{G,\mu,b,K}.
\]
The tower $(\Mscr_{G,\mu,b,K})_{K \subset G(\QQ_p)}$ is by definition the local Shimura variety attached to $(G,\mu,b)$.

Scholze and Weinstein show that these local Shimura varieties generalize the generic fibers of Rapoport-Zink spaces which are moduli spaces of $p$-divisible groups defined by Rapoport and Zink in \cite{RZ_Period}. This meets the conjectures in \cite{RV_TowardsLocal}. Their proof is based on Theorem~\ref{PDiv}.


\subsection{Integral local Shimura varieties}

Similarly, as Rapoport-Zink spaces are by definition generic fibers of formal schemes over $\Spf O_{\Ebreve}$, one also wants the local Shimura varieties constructed above to have integral models over $O_{\Ebreve}$. They should depend on the choice of a smooth model $\Gcal$ of $G$ defined over $\ZZ_p$ that is assumed to be parahoric\footnote{Or at least its connected component is assumed to be parahoric.}. Scholze and Weinstein define in \cite[XXV]{SW_Berkeley} such integral models.

Their construction depends on the existence of a good local model that models the singularities of the integral local Shimura varieties. In full generality they exist only as diamonds that are conjecturally diamonds attached to normal schemes, flat and projective over $O_{\Ebreve}$, see \cite[21.4]{SW_Berkeley} for details. In fact, Pappas and Zhu defined local models which depend on some auxiliary choices in \cite{PZ_LocalModels}. One might hope that their diamonds are the diamond local models defined by Scholze and Weinstein. This would also show that the Pappas and Zhu local models are in fact independent of the choices.

Hence in general the construction of integral models still depends on these conjectures. But in a lot of special cases these integral models are unconditionally shown to exist in \cite[25.4, 25.5]{SW_Berkeley}.


\part{Further achievements}

\section{Topological cyclic homology}\label{TCH}

Topological cyclic homology relies on the topological theory of Hochschild homology. Roughly, topological Hochschild homology (THH) is obtained by replacing the ring $\ZZ$ with the sphere spectrum $\SS$ in the definition of Hochschild homology. This yields a better integral theory. It attaches to an associative ring spectrum $A$ a spectrum ${\rm THH}(A)$ together with an $S^1$-action and an $S^1$-equivariant Frobenius map. The nature of the Frobenius map has been clarified by Nikolaus and Scholze in \cite{NS_TCH}, where they construct topological Hochschild homology as a cyclotomic spectrum. The negative (resp.~periodic) topological cyclic homology is defined as the homotopy fixed points (resp.~the Tate construction) of the $S^1$-action on ${\rm THH}(A)$. Topological cyclic homology is then given as the equalizer of the canonical map and the Frobenius from negative topological cyclic homology to $p$-completed periodic topological cyclic homology.

The relevance for algebraic geometry is furnished by the work of Bhatt-Morrow-Scholze in \cite{BMS_THH} where it is proved that B\"okstedt periodicity holds for every perfectoid ring. They define a weight filtrations of topological cyclic homology and relate the graded pieces with $p$-adic cohomology theories which now can be viewed as prismatic cohomology.


\section{Condensed mathematics}

Putting a topology onto algebraic objects usually destroys their nice algebraic properties: For example, topological abelian groups do not define an abelian category. Clausen and Scholze develop a theory that resolves this issue by replacing (and arguably improving) the notion of a topological space (\cite{CS_Condensed}). This starts with the observation that every compact Hausdorff space is the quotient of a profinite space. They define a site whose underlying category is the category of profinite spaces and whose coverings are given by finite families of arbitrary continuous maps $f_i\colon S_i \to S$ with $\bigcup_i f_i(S_i) = S$. This site is also the pro-\'etale site of a point, i.e., the pro-\'etale site of the spectrum of an algebraically closed field (Subsection~\ref{PROETSCHEME}), and hence is denoted by $*_{\rm pro\acute{e}t}$. A \emph{condensed set} is by definition a sheaf of sets on this site.

Attaching to a topological space $X$ the sheaf that sends $S$ to the set of continuous functions $S \to X$ yields a faithful functor from the category of topological spaces to the category of condensed sets. If restricted to the category of compactly generated spaces\footnote{The category of compactly generated spaces is the ``convenient category for homotopy theory'' (often one also adds the separation axiom of being weakly Hausdorff). Recall that a topological space $X$ is \emph{compactly generated} if a subset $A$ of $X$ is closed if and only if for every continuous map $f\colon S \to X$ with $S$ compact Hausdorff (or, equivalently, profinite) $f^{-1}(A)$ is closed. Every first countable space and every Hausdorff locally compact space is compactly generated.} the functor is fully faithful. Hence one may consider the category of condensed sets as an enlargement of the category of ``convenient'' topological spaces.

One obtains notions of condensed (abelian) groups or condensed rings by considering sheaves of (abelian) groups or rings. The advantage of the enlarged category of condensed abelian groups over the category of topological abelian groups is the fact that the former is automatically abelian as every category of sheaves of abelian groups on any site. This is not the case for the category of topological abelian groups in which exist bijective continuous group isomorphisms that are not homeomorphisms (e.g., the identity from $\RR$ endowed with the discrete topology to $\RR$ endowed with the standard topology).

Clausen and Scholze define the notion of completion for condensed abelian groups. This works well in the non-archimedean setting. For instance, they obtain an abelian category of certain condensed (called \emph{solid}) $\QQ_p$-vector spaces which is equipped with a tensor product and contains all Fr\'echet spaces over $\QQ_p$. The issue of generalizing topological vector spaces over $\RR$ is more complicated and has not yet appeared in written form.

By replacing for a scheme $X$ the usual derived category $D(\Oscr_X)$ of the category of $\Oscr_X$-modules by a derived category of certain condensed modules, Clausen and Scholze define for every separated finite type morphism $f$ a ``direct image with compact support'' functor $f_{!}$ and obtain a full six-functor formalism for coherent sheaves on schemes.


\end{document}